\documentclass[12pt]{article}

\usepackage{amssymb}
\usepackage{amsmath}
\usepackage{amsbsy}
\usepackage{amscd}
\usepackage{amsfonts}
\usepackage{amsthm}
\usepackage{mathrsfs}
\usepackage{verbatim}
\usepackage[colorlinks]{hyperref}
\usepackage{fullpage}
\usepackage{mathdots}
\usepackage{graphicx,subfigure}
\usepackage[english]{babel}
\usepackage[utf8]{inputenc}
\usepackage{tikz}

\usetikzlibrary{backgrounds,fit, matrix}
\usetikzlibrary{positioning}
\usetikzlibrary{calc,through,chains}
\usetikzlibrary{arrows,shapes,snakes,automata, petri}
\usepackage{authblk}

\theoremstyle{definition}
\newtheorem{Theorem}{Theorem}[section]
\newtheorem{Corollary}[Theorem]{Corollary}
\newtheorem{Lemma}[Theorem]{Lemma}
\newtheorem{Proposition}[Theorem]{Proposition}
\newtheorem{Example}[Theorem]{Example}

\newtheorem{Remark}[Theorem]{Remark}
\newtheorem{Definition}[Theorem]{Definition}
\newtheorem{Notation}[Theorem]{Notation}

\newtheorem{Conjecture}[Theorem]{Conjecture}

\newtheorem{Question}[Theorem]{Question}

\newcommand{\C}{{\mathbb C}}
\newcommand{\Z}{{\mathbb Z}}
\newcommand{\Q}{{\mathbb Q}}
\newcommand{\R}{{\mathbb R}}

\newcommand{\N}{{\mathbb N}}

\newcommand{\mc}[1]{\mathcal{#1}}

\newcommand{\mf}[1]{\mathfrak{#1}}

\newcommand{\mbf}[1]{\mathbf{ #1}}

\newcommand{\mtt}[1]{\mathtt{#1}}

\begin{document}

\title{On Generalized Wilf Conjectures}

\author[1]{Mahir Bilen Can}
\author[2]{Naufil Sakran}

\affil[1]{{\small Tulane University, New Orleans LA; \href{mailto:mahirbilencan@gmail.com}{mahirbilencan@gmail.com}}}

\affil[2]{{\small Tulane University, New Orleans LA; \href{mailto:nsakran@tulane.edu}{nsakran@tulane.edu}}}


\normalsize

\date{\today}
\maketitle

\begin{abstract}

We investigate complement-finite submonoids of the monoid of nonnegative integer points of a unipotent linear algebraic group $G$.
These monoids are in general noncommutative but they specialize to the generalized numerical monoids of Cistco et al. 
We show that every unipotent numerical monoid has a unique finite minimal generating set.
We propose a generalization of the Wilf conjecture in our setting. 
We contrast our Wilf conjecture against the Generalized Wilf Conjecture.
Then we isolate two new families of unipotent numerical monoids called the {\em thick} and the {\em thin} unipotent numerical monoids. 
We prove that our Wilf conjecture holds for every thick (commutative) unipotent numerical monoid.  
Under additional assumptions on the conductors, we prove that our Wilf conjecture holds for every thin (commutative) unipotent numerical monoid.
\medskip

\medskip

\end{abstract}

\section{Introduction}

The classification of the isomorphism classes of the subgroups of $(\Z,+)$ is a rather easy problem. 
There are essentially two isomorphism classes; one class for the trivial subgroup 0, 
and the other class for the subgroups of the form $n\Z$, where $n\in \Z\setminus \{0\}$. 
In contrast, the classification of the isomorphism classes of the submonoids of $\N$ is a deep problem. 
We pin the starting point of this classification to the following concept. 
Let $S$ be a submonoid of a monoid $M$.
We assume that $S$ is not isomorphic to $M$.
If the complement of $S$ in $M$ is a finite set, then we call $S$ a {\em complement-finite submonoid}.
For $M=\N$, the complement-finite submonoids are called the {\em numerical monoids}. 

It is peculiar to $\N$ that if an infinite submonoid $S'\subset \N$ is not a complement-finite submonoid, then $S'$ is isomorphic to a complement-finite submonoid of $\N$.
In general, for other families of infinite monoids this property does not hold.
The reader will recognize in retrospect that this is one of the starting points of our research. 
For now let us continue with $\N$ as our ambient monoid. 
It is now evident that for understanding the nature of a submonoid $S$ of $\N$,
we need to understand the complement $\N \setminus S$. 
The aim of this manuscript is to offer a generalization of a conjecture of Wilf, 
which ties together several important invariants of these  complements.
Let us introduce its ingredients: 
\begin{center}
\begin{tabular}{llll}
$\mtt{F}(S)$ &:=& $\max (\N \setminus S)$, & \qquad called the {\em Frobenius number} of $S$;\\
$\mtt{c}(S)$ &:=& $\mtt{F}(S)+1$, & \qquad called the {\em conductor} of $S$;\\
$\mtt{g}(S)$ &:=& $|\N \setminus S|$, & \qquad called the {\em genus} of $S$.
\end{tabular}
\end{center}
Let us call the elements of $S$ in $[0,\mtt{F}(S)]$ the {\em sporadic elements} of $S$. 
Then the difference $\mtt{c}(S)-\mtt{g}(S)$, denoted $\mtt{n}(S)$, is nothing but the number of sporadic elements of $S$.
We know from elementary algebra that every submonoid of $\N$ is finitely generated. 
The minimum number of generators of a numerical semigroup $S$ is called the {\em embedding dimension} of $S$, denoted $\mtt{e}(S)$.
In his short but insightful article~\cite{Wilf78}, Wilf posed a question 
that is equivalent to the following inequality: 
\begin{align*}
 \mtt{c}(S) \leq \mtt{e}(S) \mtt{n}(S),
\end{align*}
which is now commonly called the {\em Wilf conjecture}.
We will call it the {\em ordinary Wilf conjecture}.
There is a vast literature on this conjecture. 
Many papers are dedicated to the proofs of its special cases. 
However, it is still an open problem at large. 
We refer the reader to a relatively recent survey article of Delgado~\cite{Delgado} that summarizes the advances on the conjecture up to the year, 2020. 
\medskip

In this manuscript, among other things, we extend the Wilf conjecture to a setting of the group of integer points of 
a unipotent complex linear algebraic group. To spell this out more precisely, let us consider a closed subgroup $G$ of 
$GL(V)$, where $V$ is a finite dimensional vector space defined over $\C$. 
Let $1$ denote the neutral element of $G$. 
If an element $x\in G$ has the property that $x-1$ is nilpotent, then we call $x$ {\em unipotent}.
If every element of $G$ is unipotent, then $G$ is said to be a {\em unipotent linear algebraic group}. 
It is well-known that~\cite[Proposition 2.4.12]{Springer} every unipotent linear algebraic group $G$ is isomorphic to a closed 
subgroup of the unipotent group of all upper triangular unipotent $n\times n$ matrices, denoted $\mbf{U}(n,\C)$,
for some $n\in \Z_+$, where $\Z_+$ denotes the set of positive integers.
Therefore, we proceed with the assumption that $G$ is a subgroup of $\mbf{U}(n,\C)$. 
An additional assumption that we maintain throughout our manuscript is that $G$ is defined over $\Q$. 
In other words, as an algebraic variety, the defining equations of $G$ have rational coefficients.

Let $\mbf{U}(n,\Z)$ denote the group of integer points of $\mbf{U}(n,\C)$. 
Then the group of integer points of $G$, denoted $G_{\Z}$, is a subgroup of $\mbf{U}(n,\Z)$. 
The set 
\begin{align*}
\mbf{U}(n,\N):=\{ x\in \mbf{U}(n,\Z) :\ \text{all entries of $x$ are from $\N$}\}
\end{align*}
is a submonoid of $\mbf{U}(n,\Z)$. 
Likewise, the set 
\begin{align*}
G_{\N}:=\{x\in G_{\Z} : \text{all entries of $x$ are from $\N$}\}
\end{align*} 
is a submonoid of $G_{\Z}$. 
Disregarding $n$, let us call $\mbf{U}(n,\N)$ a {\em trivial numerical monoid}. 
Notice that $(\N,+)$ is a trivial numerical monoid since it is given by our simplest example, that is, $\mbf{U}(2,\N)$.
Our first nontrivial definition is the following. 

\begin{Definition}
Let $M$ be a monoid of the form $M:=G_{\N}$, where $G$ is a unipotent complex linear algebraic group.
We assume that $M$ is a finitely generated monoid.
If $S$ is a complement-finite submonoid of $M$, then we call $S$ a {\em unipotent numerical monoid in $M$}. 
\end{Definition} 

\medskip

Let $G$ denote, as above, a closed subgroup of $\mbf{U}(n,\C)$. 
We denote by $\mtt{d}_G$ the dimension of $G$ as an algebraic variety.
Let $M$ denote $G_{\N}$.
Let $S$ be a unipotent numerical monoid in $M$.
We will introduce the notions of the conductor, genus, and the sporadicity of $S$ relative to $M$.
We call these numbers, including the ``embedding dimension'' which will be introduced below, the {\em basic invariants of $S$ relative to $M$}.  
If we specialize to the case of $\mbf{U}(2,\N)$, then the basic invariants will agree with the notions that their names indicate. 
For $k\in \Z_+$, let $\mbf{U}(n,\N)_{k}$ denote the following subset of $\mbf{U}(n,\N)$: 
\begin{align*}
\mbf{U}(n,\N)_{k}:= \{ (x_{ij})_{1\leq i,j\leq n}\in \mbf{U}(n,\N):\ k\leq \max_{1\leq i,j\leq n} x_{ij}  \} \cup \{\mbf{1}_n\},
\end{align*}
where $\mbf{1}_n$ is the identity matrix of size $n$.  
We will show later that $\mbf{U}(n,\N)_{k}$ is indeed a submonoid of $\mbf{U}(n,\N)$. 
We call it the {\em $k$-th fundamental submonoid of $\mbf{U}(n,\N)$}. 
Notice that for every $k\in \Z_+$ we have $\mbf{U}(n,\N)_{k+1}\subset \mbf{U}(n,\N)_{k}$ and that $\mbf{U}(n,\N)_1 = \mbf{U}(n,\N)$. 
The {\em generating number of $S$ relative to $M$}, denoted $\mtt{r}_M(S)$, is the smallest positive integer 
$k\in \Z_+$ such that $\mbf{U}(n,\N)_{k}\cap M\subseteq S$.
We then define 
\begin{center}
\begin{tabular}{llll}
$\mtt{c}_M(S)$ &:=& $\mtt{r}_M(S)^{\mtt{d}_G}$, & \qquad called the {\em conductor of $S$ relative to $M$};\\
$\mtt{g}_M(S)$ &:=& $|M \setminus S|$, & \qquad called the {\em genus of $S$ relative to $M$};\\
$\mtt{n}_M(S)$ &:=& $|\{\mbf{1}_n\}\cup S \setminus \mbf{U}(n,\N)_{\mtt{r}_M(S)}|$, & \qquad called the {\em sporadicity of $S$ relative $M$}.
\end{tabular}
\end{center}
The elements of the set $\{\mbf{1}_n\}\cup S \setminus \mbf{U}(n,\N)_{\mtt{r}_M(S)}$ are called the {\em sporadic elements of $S$ relative to $M$}.
Next, we define the {\em embedding dimension of $S$}, denoted $\mtt{e}(S)$, as follows:
\begin{align*}
\mtt{e}(S):= \min \{ |\mc{G} | :\ \text{ $\mc{G}$ is a generating set for $S\setminus \{\mbf{1}_n\}$}\}.
\end{align*}
The definition of embedding dimension raises the question of whether there exists a unique finite generating set for $S$.
We answer this question affirmatively in our first main result. 
Before we state our first theorem, we present two simple examples demonstrating that the monoid of $\N$-valued points of an algebraic group need not be finitely generated but most notions we introduced are defined on them.

\begin{Example}\label{E:notfg}
Let $
G:= 
\left\{
\begin{bmatrix}
1 & a & b \\
0 & 1 & a \\
0 & 0 & 1
\end{bmatrix} :\
a,b\in \C
\right\}.
$
Then $G$ is a two dimensional nonabelian unipotent algebraic group. 
Its submonoid $G_\N$ is given by 
\begin{align*}
M:= G_\N=
\left\{
x_{(a,b)}:=
\begin{bmatrix}
1 & a & b \\
0 & 1 & a \\
0 & 0 & 1
\end{bmatrix} :\
a,b\in \N
\right\}.
\end{align*}
We claim that $M$ is not finitely generated. 
Indeed, for any $a\in \N$, the element $x_{(a,0)}$ cannot be written as a product of two non-identity elements in $M$. 
In other words, every generating set of $M$ contains the infinite set $\{ x_{(a,0)} :\ a\in \N\}$. 
\end{Example}

Next, we give an example of a finitely generated complement-finite unipotent submonoid.

\begin{Example}\label{E:introfg}
Let $
G:= 
\left\{
x_{(a,b,c,d,e)}:=
\begin{bmatrix}
1 & a & b & c\\
0 & 1 & 0 & d \\
0 & 0 & 1 & e \\
0 & 0 & 0 & 0 
\end{bmatrix} :\
a,b,c,d,e\in \C
\right\}.
$
Then $G$ is a five dimensional nonabelian unipotent algebraic group. 
Let $M:=G_\N$. 
Let $S$ denote the following complement-finite submonoid of $M$:
\begin{align*}
S:=  \{x_{(a,b,c,d,e)}\in M:\ (a,b,c,d,e)\neq (0,0,0,0,i) \text{ for $i\in \{1,2,3,4,5\}$} \}.
\end{align*}
Then the smallest positive integer $k\in \Z_+$ such that $\mbf{U}(4,\N)_k \cap M \subseteq S$ is $k=6$.
Hence, the generating number of $S$ relative to $M$ is given by $\mtt{r}_M(S) = 6$.
It follows that the conductor of $S$ relative to $M$ is $\mtt{c}_M(S) = 6^5 =7776$. 
Clearly, the genus of $S$ relative to $M$ is 5.
Let us calculate the sporadicity of $S$ relative to $M$: 
\begin{align*}
\mtt{n}_M(S) &= |\{\mbf{1}_4\}\cup S \setminus \mbf{U}(4,\N)_6| \\
&= | \{ x_{(a,b,c,d,e)} \in S :\ \max \{a,b,c,d,e\} <6 \} | \\
&= 6^5 - 5\\ 
&= 7771.
\end{align*}
Finally, we observe that the set 
\begin{align*}
\{ x_{(1,0,0,0,0)}, x_{(0,1,0,0,0)},x_{(0,0,1,0,0)}, x_{(0,0,0,1,0)} \} \cup \{ x_{(0,0,0,0,j)}: j\in \{6,7,\dots,11\} \}
\end{align*}
is the minimal generating set for $S$.
Hence, the embedding dimension of $S$ relative to $M$ is $\mtt{e}(S) = 10$. 

\end{Example}

\begin{Theorem}\label{T:firstmain}
Let $G$ be a unipotent linear algebraic group such that $M:=G_{\N}$ is finitely generated. 
If $S$ is a unipotent numerical monoid in $M$, then $S$ is finitely generated.
Furthermore, $S$ possesses a unique minimal set of generators. 
\end{Theorem}

We now have all the ingredients in place to state our conjecture.

\begin{Conjecture}\label{C:firstconjecture}(\textbf{Unipotent Wilf Conjecture})
Let $G$ be a unipotent linear algebraic group. Let $M=G_{\N}$. 
If $S$ is a unipotent numerical monoid in $M$, then we have 
\begin{align*}
\mtt{d}_G \mtt{c}_M(S) \leq \mtt{e}(S) \mtt{n}_M(S).
\end{align*}
\end{Conjecture}

\begin{Example}
In Example~\ref{E:introfg}, we found that 
$\mtt{d}_G =5$, $\mtt{c}_M(S) = 7776$, $\mtt{e}(S)=10$, and $\mtt{n}_M(S) = 7771$.
These numbers satisfy the inequality in Conjecture~\ref{C:firstconjecture}.
\end{Example}

In this paper, we will prove our conjecture in several important special cases. 
If $M=\mbf{U}(2,\N)$, then every complement-finite monoid $S$ in $M$ can be viewed as an ordinary numerical monoid. 
Moreover, every basic invariant of $S$ relative to $M$ becomes a basic invariant of $S$ as a numerical monoid. 
Finally, since we have $\mtt{d}_G =1$, where $G=\mbf{U}(2,\C)$,
our conjecture becomes the ordinary Wilf conjecture for numerical monoids.
\medskip

The Wilf conjecture for numerical monoids has been generalized in the papers~\cite{CDFFPU, GGMAVT},
where the main foci are on the commutative monoids only.
In a sense, our generalization of the Wilf conjecture offers further generalization. 
Let us explain this casual remark. 
We will consider the following abelian group of unipotent matrices: 
\begin{align}\label{intro:P}
\mbf{P}(n,\C):=
\left\{
\begin{bmatrix}
1 & a_1 & a_2 &  \dots & a_{n-1}  \\
0 & 1 & 0 & \dots & 0 \\
0 & 0 & 1 & \dots & 0 \\
\vdots & \vdots & \vdots & \ddots & \vdots\\
0 & 0 & 0 & \dots & 1
\end{bmatrix}
:\
\{a_1,\dots, a_{n-1}\} \subset \C
\right\} \qquad (n\geq 2).
\end{align}
As a linear algebraic group, $G:=\mbf{P}(n,\C)$ is isomorphic to the additive linear group, $(\C^{n-1},+)$.
Hence, we have $\dim G = n-1=\mtt{d}_G$. 
The monoid $M:=\mbf{P}(n,\N)$ is isomorphic to $(\N^{n-1},+)$. 
Now, for $M$, our Conjecture~\ref{C:firstconjecture} becomes the following assertion. 
\begin{Conjecture}\label{C:secondconjecture}
Let $S$ be a unipotent numerical monoid in $M:=\mbf{P}(n,\N)$. 
Then we have
\begin{align*}
(n-1) k^{n-1} \leq \mtt{e}(S) \mtt{n}_M(S).
\end{align*}
where $k$ is the generating number $\mtt{r}_M(S)$ of $S$ relative to $M$. 
\end{Conjecture}

Let $\varphi$ be the monoid isomorphism defined by 
\begin{align*}
\varphi : \mbf{P}(n,\N) &\longrightarrow \N\times \cdots \times \N \\
(a_{ij})_{1\leq i,j\leq n} &\longmapsto (a_{12},a_{13},\dots, a_{1n}).
\end{align*}
By using $\varphi$, we transfer the unipotent numerical monoids in $\mbf{P}(n,\N)$ to the complement-finite submonoids in $\N^{n-1}$. 
Although a unipotent numerical monoid in $\mbf{P}(n,\N)$ is isomorphic to its image in $\N^{n-1}$, to stress the ambient monoid,
we will follow the terminology of~\cite{FPU2016}. 
Later we will relax this aspect of our language. 
We call such a complement-finite submonoid of $\N^{n-1}$ a {\em generalized numerical monoid} 
(GMS, for short). Let $T$ be a GMS. 
The complement of $T$ in $\N^{n-1}$, that is $H(T):=\N^{n-1} \setminus T$, is called the {\em hole set} of $T$.
In~\cite{FPU2016}, Cisto, Failla, and Utano show that every GMS has a unique minimal set of generators. 
Let us denote the cardinality of a minimal system of generators by $e(T)$. 
Of course, if $\varphi(S) =T$, then we have $\mtt{e}(S)=e(T)$. 
Let $\leqslant$ denote the partial order on $\N^{n-1}$ defined by $(a_1,\dots, a_{n-1}) \leqslant (b_1,\dots, b_{n-1})$
if and only if $a_i \leq b_i$ for every $i\in \{1,\dots, n-1\}$.
Let $c(T)$ and $n(T)$ denote the following numbers:
\begin{align*}
c(T) &:=
|\{ \mbf{a} \in \N^{n-1} :\ \mbf{a} \leqslant \mbf{b} \text{ for some $\mbf{b}\in H(T)$}\}|,\\
n(T) &:=
|\{ \mbf{a} \in T :\ \mbf{a} \leqslant \mbf{b} \text{ for some $\mbf{b}\in H(T)$}\}|.
\end{align*}
In the article~\cite{CDFFPU}, the authors propose the following conjecture:
\begin{Conjecture}\label{C:generalizedconjecture}(\textbf{Generalized Wilf Conjecture})
Let $T$ be a GMS in $\N^{n-1}$. 
Then we have 
\begin{align*}
(n-1) c(T) \leq e(T) n(T).
\end{align*}
\end{Conjecture}
In their article, the authors prove Conjecture~\ref{C:generalizedconjecture} for several interesting families of generalized numerical monoids. 
Furthermore, they establish a deep connection,~\cite[Theorem 5.7]{CDFFPU}, between the zero-dimensional monomial ideals and the Generalized Wilf Conjecture.
\medskip

We will compare the Generalized Wilf Conjecture and the Unipotent Wilf Conjecture. 
It is not difficult to see that if $S$ is a unipotent numerical monoid in $M:=\mbf{P}(n,\N)$, 
then we have $c(\varphi(S))\leq \mtt{c}_M(S)$ and $n(\varphi(S))\leq \mtt{n}_M(S)$. 
In other words, we always have 
\begin{align}\label{A:becomesanequality}
e(\varphi(S))n(\varphi(S)) \leq \mtt{e}(S)\mtt{n}_M(S).
\end{align}
Let us point out a family of unipotent numerical monoids where (\ref{A:becomesanequality}) is an equality.
Let $S$ be a unipotent numerical monoid in $M$ with generating number $k$.
If the all-$(k-1)$ vector $(k-1,\dots, k-1)$ is not an element of $\varphi(S)$, then it is easily seen that 
\begin{align}
c(\varphi(S)) =  k^{n-1} = \mtt{c}_S(M)\qquad\text{and}\qquad n(\varphi(S))=\mtt{n}_M(S).
\end{align}
Thus, for such a generalized numerical monoid $S$, the validity of the Unipotent Wilf Conjecture implies the validity of the Generalized Wilf Conjecture for $\varphi(S)$.

Independently of the relationship between the two conjectures, we think that the following inequality will hold true for every 
unipotent numerical monoid in $\mbf{P}(n,\N)$.
\begin{Conjecture}\label{C:relatingconjectures}
Let $M:=\mbf{P}(n,\N)$.
Let $S$ be a unipotent numerical monoid in $M$. 
Then we have 
\begin{align}\label{A:newconjecture}
\frac{ \mtt{c}_M(S)}{c(\varphi(S))} \leq \frac{ \mtt{n}_M(S)} {n(\varphi(S))}.
\end{align}
\end{Conjecture}
\medskip

Next, we want to mention a related work of Garc\'{\i}a-Garc\'{\i}a, Mar\'{\i}n-Arag\'{o}n, and Vigneron-Tenorio without providing any details. 
In~\cite{GGMAVT}, they propose another generalization of the 
Wilf conjecture, called the Extended Wilf Conjecture, for a broad family of affine semigroups called the $\mathcal{C}$-semigroups.
It turns out that every GMS is a $\mathcal{C}$-semigroup. 
In~\cite[Proposition 6.3]{CDFFPU} Cisto et.al, shows that for a GMS, the Generalized Wilf Conjecture is a stronger assertion 
than the Extended Wilf Conjecture.

\medskip

We now go back to our discussion of the basic invariants of a unipotent numerical monoid. 
We have the following question regarding the relationship between the generating number and the genus of $S$ relative to $M$.
\begin{Question}\label{Q:thirdconjecture}
Let $S$ be a unipotent numerical monoid in $M:=G_{\N}$, where $G$ is a unipotent linear algebraic group. 
Is it true that 
\begin{align}\label{A:thirdconjecture}
\left\lfloor \frac{\mtt{r}_M(S)}{2} \right\rfloor \leq \mtt{g}_M(S) <  \mtt{c}_M(S)?
\end{align}
\end{Question}
We have a partial affirmative answer for Question~\ref{Q:thirdconjecture}.

\begin{Theorem}\label{T:thirdconjecture}
Let $M\in \{\mbf{U}(n,\N),\mbf{P}(n,\N)\}$.
If $S$ is a unipotent numerical monoid in $M$, then the inequalities in (\ref{A:thirdconjecture}) hold true.
\end{Theorem}

We are now ready to describe the structure of our paper.
Meanwhile we will mention our more specific results. 
In Section~\ref{S:Finite} we solve the finite generation problem for the unipotent numerical monoids. 
We prove our first main result Theorem~\ref{T:firstmain} in that section.
Also in this section, we determine the minimal generating set for the fundamental monoids of $\mbf{U}(n,\N)$. 
The purpose of the short Section~\ref{S:Examples} is to show that the Unipotent Wilf Conjecture holds for a fundamental monoid. 
In Section~\ref{S:Third} we prove our second main result that is Theorem~\ref{T:thirdconjecture}.
In Section~\ref{S:Coordinate} we investigate the relationship between the basic invariants of a unipotent monoid and its 
coordinate monoids. 
Here, by a coordinate monoid we mean roughly the following concept.
Let $\mbf{U}_{i,j}(n,\C)$ denote the one dimensional unipotent subgroup of $\mbf{U}(n,\C)$ that is generated by the elementary unipotent matrix $E_{i,j} = (e_{i,j})_{1\leq i , j \leq n}$ whose only nonzero non-diagonal entry is given by $e_{i,j}=1$. 
Then the $(i,j)$-th coordinate submonoid of a unipotent monoid $S$ in $M:=G_\N$ is given by $S\cap \mbf{U}_{i,j}(n,\N)$.
It turns out that the coordinate monoids can be used for providing estimates for the basic invariants of the ambient unipotent.
We do this for the genera and the conductors. 
Finally, in Section~\ref{S:Thick} we introduce two new classes of unipotent numerical monoids by using the coordinate submonoids.
We call them the {\em thick monoids} and the {\em thin monoids}. 
We prove that the Unipotent Wilf Conjecture holds for every commutative thick monoid (Theorem~\ref{T:UWCforthick}). 
It is interesting to note that our proof does not assume the validity of the ordinary Wilf conjecture.
For thin monoids, we show that there is a big subfamily in $\mbf{P}(n,\N)$ for which the Unipotent Wilf Conjecture holds true (Theorem~\ref{T:SUWCforthin}). 
In particular, we observe that, for the members of this subfamily, the Unipotent Wilf Conjecture is equivalent to the Generalized Wilf conjecture. 
Unlike our hypothesis of Theorem~\ref{T:UWCforthick}, we have the ordinary Wilf conjecture as part of our hypothesis in Theorem~\ref{T:SUWCforthin}.

\section{Finite Generation}\label{S:Finite} 

In this section we present our preliminary observations regarding the finite generation properties of the unipotent numerical monoids.
First, we setup our notation. 

\begin{Definition}
Let $P$ be a subset of a monoid $M$. 
Let $*$ denote the multiplication in $M$. 
We denote by $\langle P \rangle$ the submonoid generated by $P$ and the unit $e$ of $M$, that is, 
\begin{align*}
\langle P \rangle=\{a_1*a_2*\cdots * a_k\,:\,\{a_1,\dots, a_k\} \subset P,\ k\in \N\} \cup \{e\}.
\end{align*}
\end{Definition}

Now we will prove the first part of our first main result, Theorem~\ref{T:firstmain}.
Let us recall its statement for convenience. 
\medskip

Let $G$ be a unipotent linear algebraic group in $\mbf{U}(n,\C)$.
Let $M:=G_{\N}$.
If $S$ is a unipotent numerical monoid in $M$, then $S$ is finitely generated. 
Furthermore, $S$ possesses a unique minimal generating set.

\begin{proof}[Proof of Theorem~\ref{T:firstmain}]
We begin with showing that $S$ is a finitely generated monoid. 
The definition of a unipotent numerical monoid requires that the ambient monoid $M$ possesses a finite generating set, $A:=\{g_1,\dots, g_m\} \subset M$. Here, without loss of generality, we assume that $\mbf{1}_n \notin A$. 
Thus, every element $x\in M$ is a product of the form 
\begin{align}\label{A:multiplicand} 
x= g_{i_1}^{a_1}\cdots g_{i_l}^{a_l}
\end{align} 
for some $\{i_1,\dots, i_l\}\subseteq \{1,\dots, m\}$ and $\{a_1,\dots, a_l\}\subset \Z_+$.
Then, for each positive integer $k\in \Z_+$, let $M_k'$ denote the set of elements $x\in M$ such that a multiplicand $g_{i_s}^{a_s}$ of $x$ as in (\ref{A:multiplicand}) has exponent $a_s$ such that $k\leq a_s$.
Clearly, the union $M_k := \{ \mbf{1}_n \} \cup M_k'$ is a submonoid of $M$. 
Let us check that $M_k$ is a numerical unipotent submonoid of $M$. 
First, we observe that the elements of $M\setminus M_k$ are precisely the elements $y\in M$ such that 
$y=g_{i_1}^{a_1}\cdots g_{i_l}^{a_l}$, where $\max \{a_1,\dots, a_l\} <k$. 
This argument shows that $M_k$ is a complement-finite submonoid of $M$. 
Next, we observe that $\bigcup_{j=1}^l \{g_j^k,\dots, g_j^{2k}\}$ is a finite generating set for $M_k$. 
Hence, $M_k$ is finitely generated. 
Therefore, $M_k$ is a numerical unipotent submonoid of $M$. 

We notice that the $M_k$'s are nested, $M_1=M \supsetneq M_2 \supsetneq M_3 \supsetneq \cdots$.
Since $M\setminus S$ is a finite set, there is a smallest index $k>1$ such that $M_k\subseteq S$ and $M_{k-1}\not\subset S$.
Since both $M_k$ and $S$ are complement-finite monoids, we see that $S\setminus M_k$ is a finite set. 
It follows that the union 
\begin{align*}
\mc{G} :=(S\setminus M_k) \ \bigcup \ \left( \bigcup_{j=1}^l \{g_j^k,\dots, g_j^{2k}\} \right)
\end{align*}
is a finite generating set for $S$. 
This finishes the proof of the fact that $S$ is a unipotent numerical monoid in $M$.

We now know that $S$ possesses at least one finite generating set, $\mc{G}$, as above.
It remains to show that $\mc{G}$ contains the unique minimal generating subset of $S$. 
To this end, let us define the following subset of $S$:
\begin{align*}
\mc{T}:=\{x\in S\setminus \{\mbf{1}_n\}:\ \text{$x\neq x_1\cdot x_2\cdots x_r$ for any subset $\{x_1,\dots, x_r \} \subset S\setminus\{ x,\mbf{1}_n\}$}\}.
\end{align*}
On one hand, the finite generating set $\mc{G}$ has to contain $\mc{T}$. 
Otherwise, $\mc{G}$ cannot generate $S$. 
Indeed, if we assume that there is $x\in \mc{T}\setminus \mc{G}$, then at least one entry of $x$ is greater than $2k$. 
This implies that $x\in \langle \mc{G} \rangle \setminus \mc{G}$.
This means that $x$ can be written as a product of some elements of $\mc{G}$. Hence, we see that $x\notin \mc{T}$, which is a contradiction.
On the other hand, by using an easy induction argument, we see that every element of $\mc{G}\setminus \mc{T}$ can be written as a product of some elements of $\mc{T}$. (We demonstrate this inductive argument in the computations below.)
It follows that $\mc{T}$ is the required unique minimal finite generating set for $S$. 
This finishes the proof of our first main result. 
\end{proof}

\begin{Remark}
The anonymous referee brought to our attention that some of the conclusion of the second part of the proof of our Theorem~\ref{T:firstmain} can be derived by considering the fact that $S$ is a cancellative monoid since $S$ is a submonoid of a group. 
\end{Remark}

Consider a unipotent algebraic monoid denoted by $G$. 
We mentioned in the introduction, specifically in Example~\ref{E:notfg}, that the monoid $G_\N$ is not always finitely generated. 
Unlike the monoid $G_\N$, the group $G_\Z$ is always finitely generated. In fact, a well-known result of Borel and Harish-Chandra~\cite{BorelHarish-Chandra} asserts that if $H$ is a complex linear algebraic group defined over $\Q$, its subgroup of integer points $H_\Z$ is always finitely generated.
It would be very interesting to determine the unipotent algebraic groups $G$ such that $G_\N$ is a finitely generated monoid. 
In the rest of the paper, we will consider two important examples of such groups. 
\medskip

Let us identify for $\mbf{U}(n,\N)$ ($n\geq 2$) the unique minimal system of generators. 
Let $\mbf{M}_n$ denote the monoid of $n\times n$ matrices with entries from $\C$. 
Let $\{i,j\} \subset \{1,\dots, n\}$. 
The {\em $(i,j)$-th elementary matrix} in $\mbf{M}_n$, denoted by $E_{i,j}$, is the matrix $(e_{r,s})_{1\leq r,s \leq n}$ 
such that 
\begin{align*}
e_{r,s} :=
\begin{cases}
1 & \text{ if $r=s$};\\
1 & \text{ if $r=i,\ s=j$};\\
0 & \text{ otherwise.}
\end{cases}
\end{align*}

\begin{Lemma}\label{L:generalexpression}
Let $n\geq 2$. Let $A$ be an element of $\mbf{U}(n,\N)$.
If $A$ is given by $A:=(a_{ij})_{1\leq i,j\leq n}$, then it has a unique expression of the form 
\begin{align}\label{A:generatorsforUn}
E_{n-1,n}^{a_{n-1,n}} E_{n-2,n}^{a_{n-2,n}} E_{n-2,n-1}^{a_{n-2,n-1}} E_{n-3,n}^{a_{n-3,n}}\cdots E_{1,2}^{a_{1,2}},
\end{align}
where the ordering of the elementary matrices in (\ref{A:generatorsforUn}) is given by the reverse lexicographic ordering on the pairs of indices. 
\end{Lemma}

\begin{proof}
Let us first show the uniqueness. 
Let $A=(a_{i,j})_{1\leq i,j\leq n}$ an element of $\mbf{U}(n,\N)$.
Let us assume that $A$ has two expressions as in
\begin{align}\label{A:Aisgivenbytwo}
E_{n-1,n}^{a_{n-1,n}} E_{n-2,n}^{a_{n-2,n}} E_{n-2,n-1}^{a_{n-2,n-1}} E_{n-3,n}^{a_{n-3,n}}\cdots E_{1,2}^{a_{1,2}}
=
E_{n-1,n}^{b_{n-1,n}} E_{n-2,n}^{b_{n-2,n}} E_{n-2,n-1}^{b_{n-2,n-1}} E_{n-3,n}^{b_{n-3,n}}\cdots E_{1,2}^{b_{1,2}}.
\end{align}
The matrix $E_{n-1,n}^{a_{n-1,n}}$ is the unique factor in (\ref{A:Aisgivenbytwo}) that adds the entry $a_{n-1,n}$ to the factor 
$A':=E_{n-2,n}^{a_{n-2,n}} E_{n-2,n-1}^{a_{n-2,n-1}} E_{n-3,n}^{a_{n-3,n}}\cdots E_{1,2}^{a_{1,2}}$. 
Indeed, it is easy to check that the $(n-1,n)$-th entry of $A'$ is 0. 
Now, since $E_{n-1,n}^{a_{n-1,n}}$ is invertible, we see from the equation in (\ref{A:Aisgivenbytwo}) that $a_{n-1,n}$ must be equal to $b_{n-1,n}$. 
We proceed along the same lines to show that $a_{n-2,n}=b_{n-2,n}$. 
First we remove $E_{n-1,n}^{a_{n-1,n}}$ from both sides of the equation in (\ref{A:Aisgivenbytwo}).
Let $A''$ denote $E_{n-2,n-1}^{a_{n-2,n-1}} E_{n-3,n}^{a_{n-3,n}}\cdots E_{1,2}^{a_{1,2}}$.
It is easy to check that the $(n-2,n)$-th entry of $A''$ is 0.
Therefore, the equation 
\begin{align*}
E_{n-2,n}^{a_{n-2,n}} E_{n-2,n-1}^{a_{n-2,n-1}} E_{n-3,n}^{a_{n-3,n}}\cdots E_{1,2}^{a_{1,2}}
=
E_{n-2,n}^{b_{n-2,n}} E_{n-2,n-1}^{b_{n-2,n-1}} E_{n-3,n}^{b_{n-3,n}}\cdots E_{1,2}^{b_{1,2}}.
\end{align*}
implies that $a_{n-2,n}=b_{n-2,n}$. Continuing in this manner by the reverse lexicographic ordering on the pairs of indices,
we see that $a_{i,j}=b_{i,j}$ for every $1\leq i < j \leq n$. 
This finishes the proof of uniqueness. 
But our algorithmic proof of the uniqueness shows also that $A$ can be written as a product of the elements of the set 
$\{ E_{i,j} :\ 1\leq i < j \leq n\}$.
Hence, the proof of our lemma is finished.
\end{proof}

\begin{Proposition}\label{P:mgforUn}
The set  $\{ E_{i,j} :\ 1\leq i < j \leq n\}$ is the unique minimal generating set for $\mbf{U}(n,\N)$.
\end{Proposition}

\begin{proof}
In Lemma~\ref{L:generalexpression}, we already showed that $\{ E_{i,j} :\ 1\leq i < j \leq n\}$ is a generating set.
Let us assume towards a contradiction that an elementary matrix $E_{i,j}$ ($1\leq i<j\leq n$) can be written as a product of 
some other elements from $\mbf{U}(n,\N)$. 
By using the expressions (\ref{A:generatorsforUn}) for the factors of $E_{i,j}$, we assume that it can be written in the form 
\begin{align}\label{A:Eijintermsof}
E_{i,j} = E_{i_1,j_1}^{x_{i_1,j_1}} \cdots E_{i_r,j_r}^{x_{i_r,j_r}}
\end{align}
for some multiset of pairs of indices $\{\!\!\{ (i_1,j_1),\dots, (i_r,j_r)\}\!\!\}$ and for some multiset of exponents,
$\{\!\!\{ x_{i_1,j_1},\dots, x_{i_r,j_r}\}\!\!\}$.
Of course, by our assumption, these exponents are all positive integers. 
But in a product of the form (\ref{A:Eijintermsof}), for each factor $E_{i_s,j_s}^{x_{i_s,j_s}}$, the $(i_s,j_s)$-th entry of the product has the exponent $x_{i_s,j_s}$ as a summand. In other words, unless the multiset $\{\!\!\{ (i_1,j_1),\dots, (i_r,j_r)\}\!\!\}$ is equal to the set $\{(i,j)\}$,
the right hand side of (\ref{A:Eijintermsof}) can not be equal to the left hand side. 
This gives us the desired contradiction. 
In conclusion, we see that the basis $\{ E_{i,j} :\ 1\leq i < j \leq n\}$ is minimal. Hence, 
by Theorem~\ref{T:firstmain}, it is the unique minimal generating set for $\mbf{U}(n,\N)$.
This finishes the proof of our proposition.
\end{proof}

Our previous proposition yields many examples of unipotent numerical monoids. 
For $(i,j) \in \{1,\dots, n\}\times \{1,\dots, n\}$, let $X_{i,j}$ denote the $(i,j)$-th coordinate function on $\mbf{U}(n,\C)$ defined 
by 
\begin{align*}
X_{i,j}((a_{k,l})_{1\leq k,l\leq n}) = 
\begin{cases}
a_{i,j} & \text{ if $(i,j) = (k,l)$},\\
0 & \text{ otherwise}
\end{cases}\qquad \text{where $(a_{k,l})_{1\leq k,l\leq n}\in \mbf{U}(n,\C)$}.
\end{align*}
In this notation, we have the following consequence of Proposition~\ref{P:mgforUn}. 
\begin{Corollary}\label{C:zerocoordinates}
Let $A$ be a subset of $\{ (i,j)\in \{1,\dots, n\}\times \{1,\dots, n\} :\ 1\leq i < j \leq n\}$, where $3\leq n$, such that the set 
\begin{align*}
G := \{ (a_{k,l})_{1\leq k,l\leq n}\in \mbf{U}(n,\C) :\ X_{i,j}((a_{k,l})_{1\leq k,l\leq n})=0 \text{ for every $(i,j)\in A$}\}.
\end{align*}
is an algebraic subgroup of $\mbf{U}(n,\C)$. 
Then the monoid $G_\N$ is finitely generated. 
\end{Corollary}
\begin{proof}
The arguments of the proof of Proposition~\ref{P:mgforUn} shows also that the unique minimal generating set for $G_\N$ is given by the set $\{ E_{i,j} :\ (i,j)\notin A\}$. 
In particular, these arguments shows that $G_\N$ has a finite generating set. 
\end{proof}

\begin{Remark}
Let $G$ be a unipotent algebraic monoid as in Corollary~\ref{C:zerocoordinates}. 
There are many examples of submonoids $S\subseteq G_\N$ such that $S$ is not finitely generated. 
For example, if $G$ is equal to $\mbf{P}(3,\C)$, which is defined as in (\ref{intro:P}), then the following submonoid of $\mbf{P}(3,\N)$ is not finitely generated: 
\begin{align*}
S:=
\left\{
\begin{bmatrix}
1&a&b\\0&1&0\\0&0&1
\end{bmatrix} :\
(a,b) \in \Z_+\times \Z_+ \cup \{(0,0)\}
\right\}.
\end{align*}
To show this, we observe that $S$ is isomorphic to
\begin{align*}
S:=\{(x,y)\in \Z_+\times \Z_+ :\  xy\neq 0 \} \cup \{ (0,0)\}.
\end{align*}
Clearly, the elements $(m,1)\in S$, where $m\in \Z_+$, cannot be written as a sum of two other nonzero elements of $S$. 
Hence, $S$ does not possess any finite generating set. 

Alternatively, to reach to the same conclusion, we can use the fact that $S$ is not a complement-finite submonoid of $\mbf{P}(n,\C)$.
Therefore, by Theorem~\ref{T:firstmain} and Corollary~\ref{C:zerocoordinates}, it cannot be a finitely generated submonoid of $G_\N$. 
\end{Remark}
\medskip

We proceed to show that $\mbf{U}(n,\N)$ is filtered by a distinguished family of submonoids by using the ideas of the proof of Theorem~\ref{T:firstmain}.
\begin{Lemma}\label{L:itisamonoid}
Let $n\geq 2$. 
For $k\in \Z_+$, let $\mbf{U}(n,\N)_{k}$ denote the subset of $\mbf{U}(n,\N)$ that is defined by 
\begin{align*}
\mbf{U}(n,\N)_{k}:=\{ \mbf{1}_n\}\cup \{ (x_{ij})_{1\leq i,j\leq n} \in \mbf{U}(n,\N) :\  k\leq \max_{1\leq i,j\leq n} x_{ij} \}.
\end{align*}
Then $\mbf{U}(n,\N)_{k}$ is a unipotent numerical monoid in $\mbf{U}(n,\N)$.
\end{Lemma}

\begin{proof}
Let $A$ and $B$ be two elements of $\mbf{U}(n,\N)_{k}$. 
We will show that $AB\in \mbf{U}(n,\N)_{k}$. 
Let us assume that the entries of $A$ and $B$ are given by 
\begin{align*}
A=(a_{ij})_{1\leq i,j \leq n}\qquad\text{and}\qquad B=(b_{ij})_{1\leq i,j\leq k}.
\end{align*}
Let $C:=AB$.
Since the $(i,j)$-th entry of $C$, denoted $c_{ij}$, is defined by the formula $c_{ij}= \sum_{r=1}^n a_{ir}b_{rj}$, 
we see that $a_{ii}b_{ij} + a_{ij}b_{jj}$ is a summand of $c_{ij}$. 
In other words, we always have 
\begin{align*}
a_{ij}+b_{ij} \leq c_{ij} \qquad \text{for every $\{i,j\}\subset \{1,\dots, n\}$}.
\end{align*}
It follows that $k \leq \max C$, whence $C\in \mbf{U}(n,\N)_{k}$. 
Since we have the identity matrix is contained in $\mbf{U}(n,\N)_{k}$ as well, 
we see that $\mbf{U}(n,\N)_{k}$ is a submonoid of $\mbf{U}(n,\N)$. 
At the same time, the complement $\mbf{U}(n,\N)\setminus \mbf{U}(n,\N)_k$ consists of matrices $(x_{ij})_{1\leq i,j\leq n} \in \mbf{U}(n,\N)$ 
such that $\max_{1\leq i,j\leq n} x_{ij} <k$.
Hence, we see that $\mbf{U}(n,\N)_k$ is a complement-finite monoid. 
This finishes the proof of our assertion.
\end{proof}

\begin{Remark}
For every $n\geq 2$, the unipotent monoid $\mbf{U}(n,\N)$ is a filtered monoid. 
Indeed, for $k\in \N$, let $M_k$ denote $\mbf{U}(n,\N)_{k+1}$.
Then we have the containments, 
\begin{align}\label{A:descendingchain}
\mbf{U}(n,\N) = M_0 \supset M_1 \supset M_2 \supset \cdots.
\end{align}
\end{Remark}

Let $k\geq 2$, and define 
\begin{align*}
\mc{Q}_k:= \{ (x_{ij})_{1\leq i,j\leq n}  \in \mbf{U}(n,\N)  :\ k\leq \max_{1\leq i,j\leq n} x_{ij} < 2k \}.
\end{align*}

\begin{Lemma}\label{L:Qkisags}
The set $\mc{Q}_k$ is a generating set for $\mbf{U}(n,\N)_k$. 
\end{Lemma}

\begin{proof}
The proof follows from Lemma~\ref{L:generalexpression} and the division algorithm applied to the exponents. 
\end{proof}

We now sieve $\mc{Q}_k$ to extract a minimal generating set from it. 
\begin{align}\label{A:mingensU}
\mc{U}_k:= \{ X  \in \mc{Q}_k :\ \text{all but one entry of $X$ is not contained in $\{0,\dots, k-1\}$} \}.
\end{align}
It is easy to see that the number of elements of $\mc{U}_k$ is ${n\choose 2} k^{{n\choose 2}}$. 
Indeed, to construct an element $X\in \mc{U}_k$, we have exactly ${n\choose 2}$ coordinates to place some numbers from 
the sets $\{0,\dots, k-1\}$ and $\{k,\dots, 2k-1\}$ under the condition that, from the latter set, we are allowed to use only one element $\alpha \in\{k,\dots, 2k-1\}$.
Of course, we can choose $\alpha$ in $|\{k,\dots, 2k-1\}|$ different ways, and it can be placed in one the ${n\choose 2}$ different coordinates. 
Now there are ${n\choose 2}-1$ remaining coordinates. 
In each of these coordinates we have $k$ elements to choose from $\{0,\dots,k-1\}$.
Hence, we see that in total there ${n\choose 2} k \cdot k^{{n\choose 2}-1}$ possible ways of constructing our $X$. 
\medskip

\begin{Example}\label{E:24elements}
Let $n=3$ and $k=2$. Here are the 24 elements of $\mc{U}_2$:
\begin{center}
\begin{tabular}{llllllll}
$E_{12}^2$,& $E_{13}  E_{12}^2$, & $E_{23}  E_{12}^2$,&  $E_{23} E_{13} E_{12}^2$ , & $E_{12}^3$, & $E_{13}  E_{12}^3$,& $E_{23} E_{12}^3$, & $E_{23} E_{13}  E_{12}^3$,\\
$E_{13}^2$, & $E_{13}^2  E_{12}$, & $E_{23}  E_{13}^2$, & $E_{23} E_{13}^2 E_{12}$ ,& $E_{13}^3$, & $E_{13}^3  E_{12}$, & $E_{23} E_{13}^3$, & $E_{23} E_{13}^3  E_{12}$,\\
$E_{23}^2$, & $E_{23}^2  E_{12}$, & $E_{23}^2  E_{13}$, & $E_{23}^2 E_{13} E_{12}$ , & $E_{23}^3$, & $E_{23}^3  E_{12}$, & $E_{23}^3 E_{13}$, &  $E_{23}^3 E_{13}  E_{12}$.
\end{tabular}
\end{center}
\end{Example}

\begin{Lemma}\label{L:edimofU}
The set $\mc{U}_k$ is the unique minimal generating set for $\mbf{U}(n,\N)_k$.
In particular, the embedding dimension of $\mbf{U}(n,\N)_k$ is given by 
\begin{align*}
\mtt{e} (\mbf{U}(n,\N)_k) = {n\choose 2}k^{{n\choose 2}}.
\end{align*}
\end{Lemma}

\begin{proof}
The proof has two parts. 
First, we will show that $\mc{U}_k$ is an {\em independent set}, that is, for every $X$ and $Y$ from $\mc{U}_k$
the product $XY$ is not an element of $\mc{U}_k$. 
This will show the minimality of $\mc{U}_k$. 
Secondly, we will show that every element of $\mc{Q}_k$ can be written as a product of some elements of $\mc{U}_k$.

Now, let us assume towards a contradiction that there exist $X$ and $Y$ in $\mc{U}_k$ such that $XY\in \mc{U}_k$. 
Since the $(i,j)$-th entry of $XY$ is bigger than or equal to the sum of the $(i,j)$-th entries of $X$ and $Y$,
we see that $XY$ has at least two entries in the set $\{k,\dots, 2k-1\}$. 
This contradicts with the definition of $\mc{U}_k$. 

To prove that every element of $\mc{Q}_k$ can be written as a product of some elements of $\mc{U}_k$, we will use a similar argument. Let $X\in \mc{Q}_k$ be given by $X:=(x_{i,j})_{1\leq i,j\leq n}$. 
Let us assume that the $(i,j)$-th and the $(l,t)$-th entries of $X$ are contained in the set $\{k,\dots, 2k-1\}$. 
By using Lemma~\ref{L:generalexpression} we express $X$ in the form
\begin{align}\label{A:factorizeX}
E_{n-1,n}^{x_{n-1,n}} E_{n-2,n}^{x_{n-2,n}} E_{n-2,n-1}^{x_{n-2,n-1}} E_{n-3,n}^{x_{n-3,n}}\cdots E_{1,2}^{x_{1,2}} = X. 
\end{align}
Recall that the ordering of the elementary matrices is given by the reverse lexicographic ordering on the pairs of indices. 
Without loss of generality, we assume that $(i,j) < (l,t)$ in the reverse lexicographic ordering. 
Then we split the expression (\ref{A:factorizeX}) as follows:
\begin{align}\label{A:factorizeX2}
(E_{n-1,n}^{x_{n-1,n}} \cdots E_{i,j}^{x_{i,j}}) (E_{r,s}^{x_{r,s}}\cdots E_{1,2}^{x_{1,2}}) = X,
\end{align}
where $(r,s)$ is the predecessor of $(i,j)$ in the reverse lexicographic ordering.  
Now both of the factors, $E_{n-1,n}^{x_{n-1,n}} \cdots E_{i,j}^{x_{i,j}}$ and $E_{r,s}^{x_{r,s}}\cdots E_{1,2}^{x_{1,2}}$
are contained $\mc{Q}_k$. 
Clearly, we can repeat this factorization until each factor has only one entry that is contained in the set $\{k,\dots, 2k-1\}$. 
In other words, we can write $X$ as a product of elements from $\mc{U}_k$. 
But this observation shows also that the elements of $\mc{U}_k$ are precisely those elements of $\mbf{U}(n,\N)_k$ that cannot be expressed as a product of some other elements from $\mbf{U}(n,\N)_k$.
In other words, every element $X\in \mc{U}_k$ has to be contained in every generating set for $\mbf{U}(n,\N)_k$.
This finishes the proof of our assertion that $\mc{U}_k$ is the unique minimal generating set for $\mbf{U}(n,\N)_k$. 
We already pointed out that the number of elements of $\mc{U}_k$ is ${n\choose 2} k^{{n\choose 2}}$. 
Hence, our proof is finished. 
\end{proof}

\section{Some Special Cases}\label{S:Examples}

Let us denote the $k$-th fundamental monoid of $\mbf{P}(n,\N)$ by $\mbf{P}(n,\N)_k$. 
In other words, $\mbf{P}(n,\N)_k=\mbf{U}(n,\N)_{k}\cap \mbf{P}(n,\N)$. 
Let us summarize the basic invariants of the unipotent numerical monoids 
$S:=\mbf{U}(n,\N)_k$ in $M:=\mbf{U}(n,\N)$ and $S':=\mbf{P}(n,\N)_k$ in $M':=\mbf{P}(n,\N)$.
We begin with the former submonoid. 
Recall the minimal generating set for $\mbf{U}(n,\N)_k$:
\begin{align}\label{A:mingensofU}
\mc{U}_k 
&=\left\{X:= (x_{ij})_{1\leq i,j\leq n}  \in \mbf{U}(n,\N) : \substack{ k\leq \max_{1\leq i,j\leq n} x_{ij} < 2k,\\
\text{all but one entry of $X$ is not contained in $\{0,\dots, k-1\}$}} \right\}.
\end{align}
Then we have
\begin{center}
\begin{tabular}{llll}
$\mtt{c}_M(S)$ &:=& $k^{{n\choose 2}}$,\\ 
$\mtt{g}_M(S)$ &:=& $k^{{n\choose 2}}-1$,\\
$\mtt{n}_M(S)$ &:=& $1$,\\
$\mtt{e}(S)$ &:=& ${n\choose 2}k^{{n\choose 2}}$. 
\end{tabular}
\end{center}
Next, we have the basic invariants for $\mbf{P}(n,\N)_k$ in $\mbf{P}(n,\N)$. 
We let $\mc{P}_k$ denote the corresponding set for $\mbf{P}(n,\N)_k$, that is, 
\begin{align}\label{A:mingensP}
\mc{P}_k 
&=\left\{X:= (x_{ij})_{1\leq i,j\leq n}  \in \mbf{P}(n,\N) : \substack{ k\leq \max_{1\leq i,j\leq n} x_{ij} < 2k,\\
\text{all but one entry of $X$ is not contained in $\{0,\dots, k-1\}$}} \right\}.
\end{align}
Then we have
\begin{center}
\begin{tabular}{llll}
$\mtt{c}_{M'}(S')$ &:=& $k^{n-1}$,\\ 
$\mtt{g}_{M'}(S')$ &:=& $k^{n-1}-1$,\\
$\mtt{n}_{M'}(S')$ &:=& $1$,\\
$\mtt{e}(S')$ &:=& $(n-1)k^{n-1}$. 
\end{tabular}
\end{center}

In the rest of this section we will verify a number of special cases of the Unipotent Wilf Conjectures from the introduction.
Let us recall its claim for the unipotent numerical submonoids of $\mbf{U}(n,\N)$:
\begin{align*}
\mtt{d}_G \mtt{c}_{\mbf{U}(n,\N)}(S) \leq \mtt{e}(S) \mtt{n}_{\mbf{U}(n,\N)}(S).
\end{align*}

\begin{Example}\label{E:WilfofU(k)s}
Let $M:=\mbf{U}(n,\N)$.
For $k\geq 2$, let $S$ denote $\mbf{U}(n,\N)_k$. 
Then we have the equalities,
\begin{align*}
{n\choose 2} \mtt{c}_{M}(S) = {n\choose 2} k^{{n\choose 2}} = \mtt{e}(S) \cdot \mtt{n}_{M}(S), 
\end{align*}
confirming the Unipotent Wilf Conjecture. 
\end{Example}
\medskip

Let us recall the statement of the specialization of our Unipotent Wilf Conjecture to $\mbf{P}(n,\N)$.  
Let $S$ be a unipotent numerical submonoid of $\mbf{P}(n,\N)$. 
Then our Conjecture~\ref{C:secondconjecture} states that 
\begin{align*}
(n-1)\mtt{c}_{\mbf{P}(n,\N)}(S) \leq \mtt{e}(S) \mtt{n}_{\mbf{P}(n,\N)}(S).
\end{align*}
\medskip

\begin{Example}\label{E:WilfofP(k)s}
Let $M'=\mbf{P}(n,\N)$. 
For $k\geq 2$, let $S'$ denote $\mbf{P}(n,\N)_k$. 
In this case, we have the equalities, 
\begin{align*}
(n-1)\mtt{c}_{M'} = (n-1) k^{n-1} = \mtt{e}(S') \cdot \mtt{n}_{M'}(S').
\end{align*}
Again, the Unipotent Wilf Conjecture holds true in this case. 
\end{Example}

\medskip

Next, we have a two dimensional example. 

\begin{Example}
Let $S$ denote the unipotent numerical monoid in $\mbf{P}(3,\N)$ that is generated by the set 
$\mc{A}:=\{E_{1,3}^2, E_{1,2}^2E_{1,3}\} \cup \mc{P}_4$, where $\mc{P}_4$ is the minimal generating set for $\mbf{P}(3,\N)_4$.
Then $\mc{P}_4= \{E_{1,2}^iE_{1,3}^j, E_{1,2}^jE_{1,3}^i:\ i\in \{0,\dots,3\},\ j\in \{4,\dots,7\}\}$.
Notice that, $\mc{A}$, as a generating set, contains many redundant elements.
To show which elements are not needed, we will visualize $S$ by viewing it as a subset of $\N\times \N$ as in Figure~\ref{F:S}.
Indeed, $\mbf{P}(3,\N)$ is isomorphic to $\N\times \N$ via the isomorphism $\varphi$ that we mentioned in the introduction.
Now, the shaded boxes in Figure~\ref{F:S} represent the elements of $S$. 
The boxes with coordinates in them correspond to the elements of the generating set $\mc{A}$. 
Finally, the blue colored boxes correspond to the minimal generators of $S$. 
\begin{figure}[htp]
\begin{center}
\scalebox{.6}{
\begin{tikzpicture}[scale=1.2]

\draw (0,0)  node[fill=olive!45!, minimum size=1cm,draw] {};

\draw (0,2)  node[fill=blue!45!, minimum size=1cm,draw] {};
\node at (0,2) {$(0,2)$};

\draw (2,1)  node[fill=blue!45!, minimum size=1cm,draw] {};
\node at (2,1) {$(2,1)$};

\draw (2,3)  node[fill=olive!45!, minimum size=1cm,draw] {};

\foreach \i in {4,...,8} {\draw (\i,0)  node[fill=olive!45!, minimum size=1cm,draw] {};};
\foreach \i in {4,...,7} {\draw (\i,0)  node[fill=blue!45!, minimum size=1cm,draw] {};};
\foreach \i in {4,...,8} {\draw (\i,1)  node[fill=olive!45!, minimum size=1cm,draw] {};};
\foreach \i in {4,...,5} {\draw (\i,1)  node[fill=blue!45!, minimum size=1cm,draw] {};};
\foreach \i in {4,...,8} {\draw (\i,2)  node[fill=olive!45!, minimum size=1cm,draw] {};};
\foreach \i in {4,...,8} {\draw (\i,3)  node[fill=olive!45!, minimum size=1cm,draw] {};};
\foreach \i in {0,...,8} {\draw (\i,4)  node[fill=olive!45!, minimum size=1cm,draw] {};};
\foreach \i in {0,...,8} {\draw (\i,5)  node[fill=olive!45!, minimum size=1cm,draw] {};};
\foreach \i in {0,...,8} {\draw (\i,6)  node[fill=olive!45!, minimum size=1cm,draw] {};};
\foreach \i in {0,...,8} {\draw (\i,7)  node[fill=olive!45!, minimum size=1cm,draw] {};};
\foreach \i in {0,...,8} {\draw (\i,8)  node[fill=olive!45!, minimum size=1cm,draw] {};};

\draw (0,5)  node[fill=blue!45!, minimum size=1cm,draw] {};
\draw (1,4)  node[fill=blue!45!, minimum size=1cm,draw] {};
\draw (1,5)  node[fill=blue!45!, minimum size=1cm,draw] {};
\draw (2,4)  node[fill=blue!45!, minimum size=1cm,draw] {};
\draw (3,4)  node[fill=blue!45!, minimum size=1cm,draw] {};

\node at (0,4) {$(0,4)$};
\node at (0,5) {$(0,5)$};
\node at (0,6) {$(0,6)$};
\node at (0,7) {$(0,7)$};

\node at (1,4) {$(1,4)$};
\node at (1,5) {$(1,5)$};
\node at (1,6) {$(1,6)$};
\node at (1,7) {$(1,7)$};

\node at (2,4) {$(2,4)$};
\node at (2,5) {$(2,5)$};
\node at (2,6) {$(2,6)$};
\node at (2,7) {$(2,7)$};

\node at (3,4) {$(3,4)$};
\node at (3,5) {$(3,5)$};
\node at (3,6) {$(3,6)$};
\node at (3,7) {$(3,7)$};

\node at (4,0) {$(4,0)$};
\node at (5,0) {$(5,0)$};
\node at (6,0) {$(6,0)$};
\node at (7,0) {$(7,0)$};

\node at (4,1) {$(4,1)$};
\node at (5,1) {$(5,1)$};
\node at (6,1) {$(6,1)$};
\node at (7,1) {$(7,1)$};

\node at (4,2) {$(4,2)$};
\node at (5,2) {$(5,2)$};
\node at (6,2) {$(6,2)$};
\node at (7,2) {$(7,2)$};

\node at (4,3) {$(4,3)$};
\node at (5,3) {$(5,3)$};
\node at (6,3) {$(6,3)$};
\node at (7,3) {$(7,3)$};

\foreach \i in {0,...,9} {\draw [dashed] (\i-.5,-.5) -- (\i-.5,9);};
\foreach \i in {0,...,9} {\draw [dashed] (-.5,\i-0.5) -- (9,\i-0.5);};

\end{tikzpicture}
}
\end{center}
\caption{A unipotent numerical monoid whose minimal generators are highlighted.}
\label{F:S}
\end{figure}
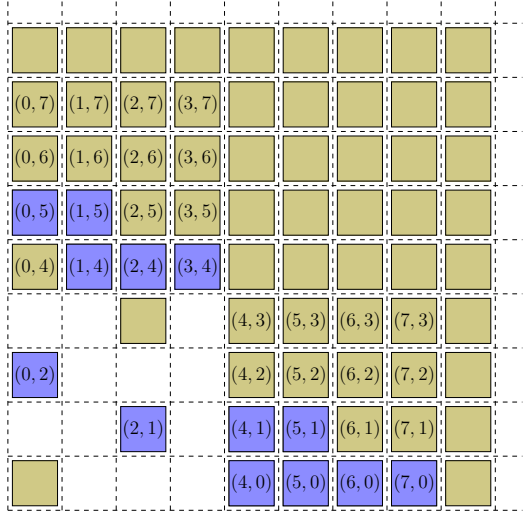
It is now easily seen that $\mtt{e}(S) = 13$, $\mtt{n}_{\mbf{P}(3,\N)}(S) = 4$, and $\mtt{c}_{\mbf{P}(3,\N)}(S)=16$. 
Then we have 
\begin{align*}
2\cdot \mtt{c}_{\mbf{P}(3,\N)}(S)  = 32\ \leq \ \mtt{e}(S) \mtt{n}_{\mbf{P}(3,\N)}(S)=52.
\end{align*}
This confirms the Unipotent Wilf Conjecture in this case.
\end{Example}

\section{Proof of Theorem~\ref{T:thirdconjecture}}\label{S:Third}

In this section, we give a partial affirmative answer to Question~\ref{Q:thirdconjecture}.
More precisely, we will prove the following statement:  
\medskip

If $S$ is a unipotent numerical monoid in $M$, where $M\in \{\mbf{U}(n,\N), \mbf{P}(n,\N)\}$, 
then the inequalities 
\begin{align*}
\left\lfloor \frac{\mtt{r}_M(S)}{2} \right\rfloor \leq \mtt{g}_M(S) <  \mtt{c}_M(S)
\end{align*}
hold true.
\medskip

In the proof of this statement, our implicit assumption will be that the generating number of $S$ relative to $M$ is at least two. 
In fact, if $\mtt{r}_M(S)=1$, then we see that $S=M$. In this case, we have 
\begin{align*}
0 = \left\lfloor \frac{\mtt{r}_M(S)}{2} \right\rfloor = \mtt{g}_M(S) <  \mtt{c}_M(S) =1.
\end{align*}
We are now ready to start our proof.

\begin{proof}[Proof of Theorem~\ref{T:thirdconjecture}.]

We will present the proof for $M=\mbf{U}(n,\N)$ only.
The case of $\mbf{P}(n,\N)$ is similar. 

We begin with introducing some new notation and terminology. 
First, we will abbreviate $\mbf{U}(n,\N)_k$ to $U(k)$. 
Now recall that $\mtt{c}_M(S)=k^{{n\choose 2}}$, where $k$ is the generating number, $\mtt{r}_M(S)$.
We interpret $k^{{n\choose 2}}$ as the number of integer points of the hypercube 
\begin{align*}
[0,k-1]\times \cdots \times [0,k-1]\subset \R^{{n\choose 2}}.
\end{align*}
This interpretation is justified by the fact that $M$ can be identified as a set (not as a monoid) with $\N^{{n\choose 2}}$.
Here, each coordinate of $\N^{{n\choose 2}}$ corresponds to the submonoid of $\mbf{U}(n,\N)$ that is generated by the elementary matrix $E_{i,j}$ ($1\leq i < j \leq n$).
Then the complement $M\setminus U(k)$ corresponds to the set of integer points of $[0,k-1]^{{n\choose 2}}$. 
Let us call this cube the {\em enveloping hypercube of $M\setminus S$}. 
Thus, we stress the fact that $\mtt{c}_M(S)$ is the volume of the enveloping hypercube.
Equivalently, we have $\mtt{c}_M(S) = |M\setminus U(k)|$.
Since the genus of $S$ relative to $M$ is the cardinality of the hole set, that is, $\mtt{g}_M(S)= |M\setminus S|$, 
the validity of the inequality $\mtt{g}_M(S) <  \mtt{c}_M(S)$ is obvious. 
\medskip

We proceed to the proof of the inequality, $\left\lfloor \frac{k}{2} \right\rfloor \leq \mtt{g}_M(S)$.
Towards a contradiction let us assume that $\left\lfloor \frac{k}{2} \right\rfloor  > \mtt{g}_M(S)$.
Here, we have a simple but crucial observation:
\medskip 

Since $k$ is the generating number of $S$ relative to $M$, the monoid $S$ does not contain the generating set $\mc{U}_{k-1}$ of $U(k-1)$. 
Hence, there exists a matrix of the form 
\begin{align*}
V := E_{n-1,n}^{v_{n-1,n}}E_{n-2,n}^{v_{n-2,n}}\cdots E_{1,2}^{v_{1,2}}  \in \mc{U}_{k-1}\setminus S,
\end{align*}
where $\{v_{i,j} :\ 1\leq i < j \leq n\}\subset \{0,\dots, k-1\}$.
Let $a:=v_{r,s}$ be the exponent such that 
\begin{align*}
\max \{ v_{1,2},\dots, v_{n-1,n} \} = a \qquad\text{and}\qquad k-1\leq a <2(k-1).
\end{align*}
All of the other elements of $\{v_{1,2},\dots, v_{n-1,n}\}$ are contained in $\{0,\dots, k-2\}$. 
But since $V\notin \mc{U}_k$, we must have $a\notin \{k,\dots, 2k-1\}$.
Hence, we conclude that $a= k-1$. 
We now have two distinct possibilities here:
\begin{enumerate}
\item[(1)] $\{v_{1,2},\dots, v_{n-1,n}\}\setminus \{a\} \neq \{0\}$, or, 
\item[(2)] $\{v_{1,2},\dots, v_{n-1,n}\}\setminus \{a\} = \{0\}$.
\end{enumerate}

Let us proceed with (1).
For this part of the proof we will use a counting argument. 
We factorize $V$ as follows. 
\begin{align*}
V := (E_{n-1,n}^{v_{n-1,n}}\cdots E_{r,s}^{a})\cdot (E_{r',s'}^{v_{r',s'}}\cdots E_{1,2}^{v_{1,2}}),
\end{align*}
where $(r',s')$ is the predecessor of $(r,s)$ in the reverse lexicographic ordering on $\N^2$. 
Since $V\notin S$, if $E_{n-1,n}^{v_{n-1,n}}\cdots E_{r,s}^{a}$ is an element of $S$, then $E_{r',s'}^{v_{r',s'}}\cdots E_{1,2}^{v_{1,2}}$ is not an element of $S$. 
Conversely, if $E_{r',s'}^{v_{r',s'}}\cdots E_{1,2}^{v_{1,2}}$ is an element of $S$, then $E_{n-1,n}^{v_{n-1,n}}\cdots E_{r,s}^{a}$ is not an element of $S$.
Another possibility here is that neither $E_{n-1,n}^{v_{n-1,n}}\cdots E_{r,s}^{a}$ nor $E_{r',s'}^{v_{r',s'}}\cdots E_{1,2}^{v_{1,2}}$ is an element of $S$. 
Likewise, if $E_{n-1,n}^{v_{n-1,n}}\cdots E_{r,s}^{a-1}$ is an element of $S$, then 
$E_{r,s}E_{r',s'}^{v_{r',s'}}\cdots E_{1,2}^{v_{1,2}}$ is not an element of $S$. 
Conversely, if $E_{r,s}E_{r',s'}^{v_{r',s'}}\cdots E_{1,2}^{v_{1,2}}$ is an element of $S$, 
then $E_{n-1,n}^{v_{n-1,n}}\cdots E_{r,s}^{a-1}$ is not an element of $S$. 
Another possibility is that neither $E_{n-1,n}^{v_{n-1,n}}\cdots E_{r,s}^{a-1}$ nor $E_{r,s}E_{r',s'}^{v_{r',s'}}\cdots E_{1,2}^{v_{1,2}}$ is an element of $S$. 
We continue in this manner for each matrix of the form 
$E_{n-1,n}^{v_{n-1,n}}\cdots E_{r,s}^{a-j}$, where $j\in \{0,\dots, k-1\}$.
At the end, we see that there are at least $k$ elements in the complement $M\setminus S$. 
This means that the genus $\mtt{g}_M(S)$ is at least $k$, which contradicts our initial assumption that 
$\left\lfloor \frac{k}{2} \right\rfloor  > \mtt{g}_M(S)$.

We now proceed with the case of (2). We will apply a similar counting argument.
In this case, our matrices are given by the powers of $E_{r,s}$ only.
Similarly to the previous case, we observe that, for $j\in \{0,\dots, k-1\}$, 
either $E_{r,s}^j$ or $E_{r,s}^{a-j}$ is not an element of $S$. 
Thus there are at least $\left\lfloor \frac{k}{2} \right\rfloor$ elements in $M\setminus S$.
Once again this contradicts with our initial assumption that $\left\lfloor \frac{k}{2} \right\rfloor  > \mtt{g}_M(S)$. This finishes the proof of our second inequality. 
Hence, the proof is complete. 
\end{proof}

\section{Coordinate Monoids}\label{S:Coordinate}

Let $n\geq 2$. 
Let $M_{i,j}$ ($1\leq i<j\leq n$) denote the monoid generated by the elementary matrix $E_{i,j}$, that is, 
\begin{align*}
M_{i,j} := \langle E_{i,j} \rangle = \{ E_{i,j}^s :\ s\in \N \} \cong \N.
\end{align*}

\begin{Definition}
Let $S$ be a submonoid of $\mbf{U}(n,\N)$. 
The {\em $(i,j)$-th coordinate submonoid of $S$}, denoted by $S_{i,j}$, is the monoid defined by 
\begin{align*}
S_{i,j} := S\cap M_{i,j}.
\end{align*}
\end{Definition}

\begin{Lemma}
Let $M$ be an element of $\{ \mbf{U}(n,\N),\mbf{P}(n,\N)\}$. 
If $S$ is a unipotent numerical monoid in $M$, then $S_{i,j}$ is a numerical monoid.
In other words, $S_{i,j}$ is isomorphic to a complement-finite submonid of $\N$. 
Here, it is understood that if $M=\mbf{P}(n,\N)$, then $i=1$.
\end{Lemma}

\begin{proof}
Since $M\setminus S$ is a finite set, the intersection $(M\setminus S)\cap M_{i,j}$ is a finite set as well.
But we have $(M\setminus S)\cap M_{i,j} = M\cap M_{i,j} \setminus S\cap M_{i,j} = M_{i,j} \setminus S_{i,j}$.
Hence, $S_{i,j}$ is a complement-finite submonoid of $M_{i,j}$. This finishes the proof of our assertion.
\end{proof}

The purpose of this section is to discuss the relationship between the basic invariants of $S$ in relation with the basic invariants of the numerical monoids $S_{i,j}$'s. 
Although our results hold for both of the ambient monoids $\mbf{U}(n,\N)$ and $\mbf{P}(n,\N)$,
we will state our theorems for the latter monoid only. 
\begin{Notation}
In the rest of this section, $M$ will stand for the commutative monoid $\mbf{P}(n,\N)$.
If the index $n$ is evident from the context, then we write $P(k)$ instead of $\mbf{P}(n,\N)_k$. 
Also, by using the isomorphism $\varphi$ from the introduction, we will view $M$ as the monoid $\N^{n-1}$ without mentioning it again. 
If $S$ is a submonoid of $M$, then we will list the coordinate submonoids of $S$ as $S_1,\dots, S_{n-1}$. 
More precisely, for $j\in \{1,\dots, n-1\}$, we set  
\begin{align*}
S_j:= S\cap (\{0\}\times \cdots \times \N\times \cdots \times \{0\}),
\end{align*}
where $\N$ appears at the $j$-th coordinate of the product.
\end{Notation}

The question that we want to answer is the following: 
\medskip
How do the basic invariants of $S$ relate to the basic invariants of the numerical monoids,
$S_1,\dots, S_{n-1}$?
\medskip

Hereafter, to make sure that $S\neq M$, we assume that the generating number of $S$ relative to $M$ is at least two.

\subsection{Conductors.}\label{SS:Conductors}

Let $k$ denote the generating number of $S$ relative to $M$. 
Then the conductor of $S$ relative to $M$ is $(n-1)k^{n-1}$. 
It is easy to check that the conductor of every coordinate submonoid $S_j$ is at most $k$. 
Thus, unless $n=2$, we have an absolute inequality,
\begin{align}\label{A:conductors}
\mtt{c}(S_1) \cdots \mtt{c}(S_{n-1}) \leq k^{n-1} < (n-1)k^{n-1} = \mtt{d}_G \mtt{c}_M(S),
\end{align}
where $G=\mbf{P}(n,\C)$.
The following inequality always holds as well: 
\begin{align}\label{A:conductors2}
\mtt{c}(S_1) +\cdots + \mtt{c}(S_{n-1}) \leq (n-1)k < (n-1)k^{n-1} = \mtt{d}_G\mtt{c}_M(S).
\end{align}
However, we want to emphasize the fact, even if the generating number of $S$ relative to $M$ is fixed, 
the conductors $\mtt{c}(S_1),\dots, \mtt{c}(S_{n-1})$ usually vary in a wide range of numbers. 
We will justify our remark by presenting several examples. 

\begin{Example}
Let $S$ denote the following submonoid of $\mbf{P}(3,\N)$:
\begin{align*}
S := \langle (1,0), (1,1),(1,2),A :\ A\in P(3) \rangle.
\end{align*}
Then we have $\mtt{r}_M(S)=3$, $\mtt{c}(S_1) =1$ and $\mtt{c}(S_2) = 3$.
\end{Example}

\begin{Example}
Let $S$ denote the monoid $P(3)\subset \mbf{P}(3,\N)$. 
Then we have $\mtt{r}_M(S)=3$, $\mtt{c}(S_1) =3$ and $\mtt{c}(S_2) = 3$.
\end{Example}

\begin{Example}
Let $S$ denote the monoid
\begin{align*}
S := \langle (2,0), (0,2),A :\ A\in P(3) \rangle \subset \mbf{P}(3,\N).
\end{align*}
Then we have $\mtt{r}_M(S)=3$, $\mtt{c}(S_1) =2$ and $\mtt{c}(S_2) =2$.
\end{Example}

These examples show that the sign of the difference 
\[
\max \{ \mtt{c}(S_1), \dots, \mtt{c}(S_{n-1}) \} - \mtt{r}_M(S)
\] 
varies in the set $\{-1,0\}$. 
In the next section we will discuss some important examples of unipotent numerical semigroups such that $\max \{ \mtt{c}(S_1), \dots, \mtt{c}(S_{n-1}) \} = \mtt{r}_M(S)$.

\subsection{Genera.}

In Theorem~\ref{T:thirdconjecture} we proved the inequalities, $\left\lfloor \frac{\mtt{r}_M(S)}{2} \right\rfloor \leq \mtt{g}_M(S) <  \mtt{c}_M(S)$.
In (\ref{A:conductors}) we found the most straightforward relationship between the conductors of 
$S$ and the $S_j$'s.
It turns out that the relationship between the genera is more complicated but leads to an interesting bound.
\medskip

Let $j\in \{1,\dots, n-1\}$.
It follows directly from the definition of $S_j$ that the hole set $M\setminus S$ contains more elements 
than the hole set $M_j \setminus S_j$. 
In other words, we have $\mtt{g}(S_j)\leq \mtt{g}_M(S)$.
In fact, similarly to (\ref{A:conductors2}), the following inequality always holds is easily seen: 
\begin{align}\label{A:additivelowerbound}
\mtt{g}(S_1)+\cdots +\mtt{g}(S_{n-1}) \leq \mtt{g}_M(S).
\end{align}
But as we will show in our next examples, unlike (\ref{A:conductors}), 
the product $\mtt{g}(S_1)\cdots \mtt{g}(S_{n-1})$ is not necessarily smaller than $\mtt{g}_M(S)$. 

\begin{Example}

In Figure~\ref{F:S2} we have two examples of unipotent numerical monoids in $\mbf{P}(3,\N)$, 
each represented by the shaded boxes in the corresponding grid. 

The one on the left is given by
\begin{align*}
S^{(1)} = \langle (1,1),(1,2),(1,3),(2,1),A :\ A\in P(4) \rangle.
\end{align*}
The one on the right, denoted $S^{(2)}$, is the fourth fundamental monoid in $\mbf{P}(3,\N)$.

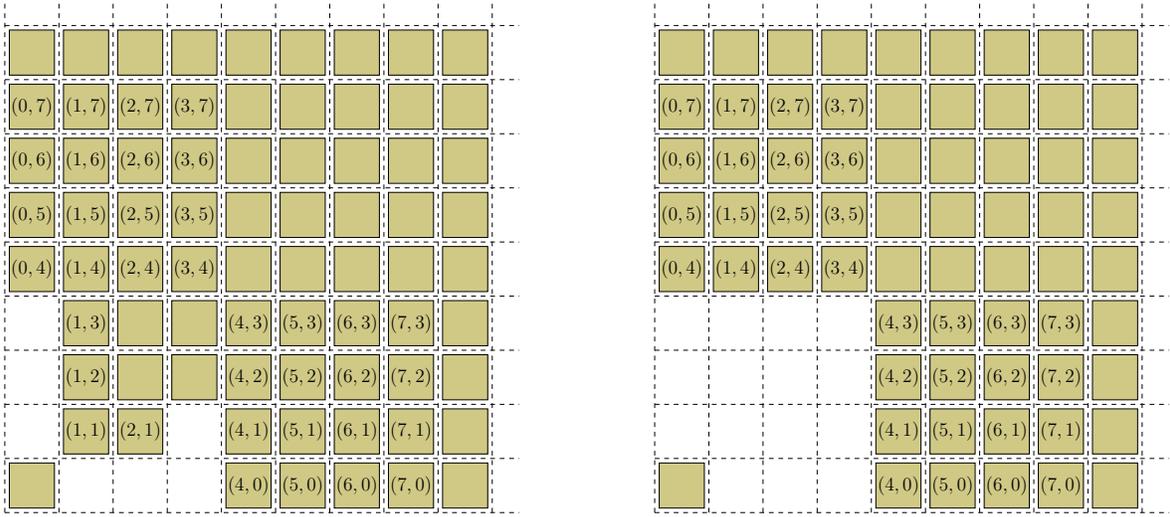
\begin{figure}[htp]
\begin{center}
\scalebox{.6}{
\begin{tikzpicture}[scale=1.2]

\begin{scope}[xshift=-6cm]
\draw (0,0)  node[fill=olive!45!, minimum size=1cm,draw] {};
\draw (1,1)  node[fill=olive!45!, minimum size=1cm,draw] {};
\draw (2,1)  node[fill=olive!45!, minimum size=1cm,draw] {};

\foreach \i in {4,...,8} {\draw (\i,0)  node[fill=olive!45!, minimum size=1cm,draw] {};};
\foreach \i in {4,...,8} {\draw (\i,1)  node[fill=olive!45!, minimum size=1cm,draw] {};};
\foreach \i in {1,...,8} {\draw (\i,2)  node[fill=olive!45!, minimum size=1cm,draw] {};};
\foreach \i in {1,...,8} {\draw (\i,3)  node[fill=olive!45!, minimum size=1cm,draw] {};};
\foreach \i in {0,...,8} {\draw (\i,4)  node[fill=olive!45!, minimum size=1cm,draw] {};};
\foreach \i in {0,...,8} {\draw (\i,5)  node[fill=olive!45!, minimum size=1cm,draw] {};};
\foreach \i in {0,...,8} {\draw (\i,6)  node[fill=olive!45!, minimum size=1cm,draw] {};};
\foreach \i in {0,...,8} {\draw (\i,7)  node[fill=olive!45!, minimum size=1cm,draw] {};};
\foreach \i in {0,...,8} {\draw (\i,8)  node[fill=olive!45!, minimum size=1cm,draw] {};};

\node at (1,1) {$(1,1)$};
\node at (1,2) {$(1,2)$};
\node at (1,3) {$(1,3)$};
\node at (2,1) {$(2,1)$};

\node at (0,4) {$(0,4)$};
\node at (0,5) {$(0,5)$};
\node at (0,6) {$(0,6)$};
\node at (0,7) {$(0,7)$};

\node at (1,4) {$(1,4)$};
\node at (1,5) {$(1,5)$};
\node at (1,6) {$(1,6)$};
\node at (1,7) {$(1,7)$};

\node at (2,4) {$(2,4)$};
\node at (2,5) {$(2,5)$};
\node at (2,6) {$(2,6)$};
\node at (2,7) {$(2,7)$};

\node at (3,4) {$(3,4)$};
\node at (3,5) {$(3,5)$};
\node at (3,6) {$(3,6)$};
\node at (3,7) {$(3,7)$};

\node at (4,0) {$(4,0)$};
\node at (5,0) {$(5,0)$};
\node at (6,0) {$(6,0)$};
\node at (7,0) {$(7,0)$};

\node at (4,1) {$(4,1)$};
\node at (5,1) {$(5,1)$};
\node at (6,1) {$(6,1)$};
\node at (7,1) {$(7,1)$};

\node at (4,2) {$(4,2)$};
\node at (5,2) {$(5,2)$};
\node at (6,2) {$(6,2)$};
\node at (7,2) {$(7,2)$};

\node at (4,3) {$(4,3)$};
\node at (5,3) {$(5,3)$};
\node at (6,3) {$(6,3)$};
\node at (7,3) {$(7,3)$};

\foreach \i in {0,...,9} {\draw [dashed] (\i-.5,-.5) -- (\i-.5,9);};
\foreach \i in {0,...,9} {\draw [dashed] (-.5,\i-0.5) -- (9,\i-0.5);};
\end{scope}

\begin{scope}[xshift=6cm]
\draw (0,0)  node[fill=olive!45!, minimum size=1cm,draw] {};
\foreach \i in {4,...,8} {\draw (\i,0)  node[fill=olive!45!, minimum size=1cm,draw] {};};
\foreach \i in {4,...,8} {\draw (\i,1)  node[fill=olive!45!, minimum size=1cm,draw] {};};
\foreach \i in {4,...,8} {\draw (\i,2)  node[fill=olive!45!, minimum size=1cm,draw] {};};
\foreach \i in {4,...,8} {\draw (\i,3)  node[fill=olive!45!, minimum size=1cm,draw] {};};
\foreach \i in {0,...,8} {\draw (\i,4)  node[fill=olive!45!, minimum size=1cm,draw] {};};
\foreach \i in {0,...,8} {\draw (\i,5)  node[fill=olive!45!, minimum size=1cm,draw] {};};
\foreach \i in {0,...,8} {\draw (\i,6)  node[fill=olive!45!, minimum size=1cm,draw] {};};
\foreach \i in {0,...,8} {\draw (\i,7)  node[fill=olive!45!, minimum size=1cm,draw] {};};
\foreach \i in {0,...,8} {\draw (\i,8)  node[fill=olive!45!, minimum size=1cm,draw] {};};

\node at (0,4) {$(0,4)$};
\node at (0,5) {$(0,5)$};
\node at (0,6) {$(0,6)$};
\node at (0,7) {$(0,7)$};

\node at (1,4) {$(1,4)$};
\node at (1,5) {$(1,5)$};
\node at (1,6) {$(1,6)$};
\node at (1,7) {$(1,7)$};

\node at (2,4) {$(2,4)$};
\node at (2,5) {$(2,5)$};
\node at (2,6) {$(2,6)$};
\node at (2,7) {$(2,7)$};

\node at (3,4) {$(3,4)$};
\node at (3,5) {$(3,5)$};
\node at (3,6) {$(3,6)$};
\node at (3,7) {$(3,7)$};

\node at (4,0) {$(4,0)$};
\node at (5,0) {$(5,0)$};
\node at (6,0) {$(6,0)$};
\node at (7,0) {$(7,0)$};

\node at (4,1) {$(4,1)$};
\node at (5,1) {$(5,1)$};
\node at (6,1) {$(6,1)$};
\node at (7,1) {$(7,1)$};

\node at (4,2) {$(4,2)$};
\node at (5,2) {$(5,2)$};
\node at (6,2) {$(6,2)$};
\node at (7,2) {$(7,2)$};

\node at (4,3) {$(4,3)$};
\node at (5,3) {$(5,3)$};
\node at (6,3) {$(6,3)$};
\node at (7,3) {$(7,3)$};

\foreach \i in {0,...,9} {\draw [dashed] (\i-.5,-.5) -- (\i-.5,9);};
\foreach \i in {0,...,9} {\draw [dashed] (-.5,\i-0.5) -- (9,\i-0.5);};

\end{scope}

\end{tikzpicture}
}
\end{center}
\caption{Two examples with different genera.}
\label{F:S2}
\end{figure}

Then the genus of $S^{(1)}$ relative to $\mbf{P}(3,\N)$ is $7$. 
Likewise, the genus of $S^{(2)}$ relative to $\mbf{P}(3,\N)$ is $15$. 
At the same time, we have 
\begin{align*}
S^{(1)}_1 = \langle (4,0)\rangle \ \cong \ S^{(1)}_2 = \langle (0,4) \rangle,
\end{align*}
hence, $\mtt{g}(S^{(1)}_1) =\mtt{g}(S^{(1)}_2) = 3$. 
It follows that 
\begin{align*}
\mtt{g}(S^{(1)}_1)\mtt{g}(S^{(1)}_2) =9 > \mtt{g}_{\mbf{P}(3,\N)}  (S^{(1)})=7.
\end{align*}

On the other hand, the genus of $S^{(2)}$ relative to $\mbf{P}(3,\N)$ is $15$. 
Also, we have $S^{(2)}_1\cong S^{(2)}_2 \cong \langle 4 \rangle$ with 
$\mtt{g}(S^{(2)}_1) = 3$. 
Thus, in this case, we have the opposite inequality, 
\begin{align*}
\mtt{g}(S^{(2)}_1) \mtt{g}(S^{(2)}_2) =9 < \mtt{g}_{\mbf{P}(3,\N)} (S^{(2)}) =15.
\end{align*}

\end{Example}

To find a good estimate for the genera, the following lemma on the sporadicity will be helpful.

\begin{Lemma}\label{L:sporadicelements}
For every unipotent numerical monoid $S$ in $M$, we have 
\begin{align*}
\mtt{n}(S_1)\cdots \mtt{n}(S_{n-1}) \leq \mtt{n}_M(S).
\end{align*}
\end{Lemma}

\begin{proof}
Let $k$ denote the generating number of $S$ relative to $M$.
Recall that the enveloping hypercube of $M\setminus S$ is the complement $M\setminus P(k)$.
Clearly, for every $j\in \{1,\dots, n-1\}$, 
the set of sporadic elements of the coordinate monoid $S_j$ is contained in the enveloping hypercube of $S$. 
Since $S$ is a monoid, all sums of the form 
\begin{align*}
a+ b\qquad\text{where $a\in S_i$ and $b\in S_j$}
\end{align*}
for $1\leq i \neq j \leq n-1$ are contained in the enveloping hypercube of $M\setminus S$.
In other words, the set of sporadic elements of $S$ relative to $M$ contains all
possible sums (with no repetitions) of the sporadic elements of $S_1,\dots, S_j$. 
This finishes the proof of our lemma.
\end{proof}

Let $\{x_1,\dots, x_{n-1}\}$ be a set of $n-1$ algebraically independent variables. 
The {\em $j$-th elementary symmetric polynomial}, for $j\in \{1,\dots, n-1\}$, is the polynomial 
$\mf{e}_j$ defined by the formula 
\begin{align*}
\mf{e}_j (x_1,\dots, x_{n-1}) = \sum_{ 1\leq i_1 < \cdots < i_j \leq {n-1}} x_{i_1}\cdots x_{i_j}.
\end{align*}

\begin{Proposition}\label{P:generacomparison}
Let $\mtt{g}_1,\dots, \mtt{g}_{n-1}$ denote the genera of $S_1,\dots,S_{n-1}$, respectively. 
Then we have 
\begin{align}\label{A:genera}
\mf{e}_1(\mtt{g}_1,\dots, \mtt{g}_{n-1})\ \leq \ \mtt{g}_M(S)\ \leq \ \sum_{j=1}^{n-1} (-1)^{j-1} k^{n-1-j}  \mf{e}_j(\mtt{g}_1,\dots, \mtt{g}_{n-1}).
\end{align}
\end{Proposition}

\begin{proof}

Since $\mf{e}_1(\mtt{g}_1,\dots, \mtt{g}_{n-1})=\sum_{j=1}^{n-1} \mtt{g}_j$, the first inequality is a restatement of (\ref{A:additivelowerbound}).
We proceed to show that the upper bound holds.

Similarly to the ordinary numerical monoids, 
for every unipotent numerical monoid $S$ in $M=\mbf{P}(n,\N)$, the sum of genus and the sporadicity gives the conductor:
\begin{align}\label{A:alwaystrue}
\mtt{n}_M(S) + \mtt{g}_M(S) = \mtt{c}_M(S) = k^{n-1},
\end{align}
where, $k$ denotes, as before, the generating number of $S$ relative to $M$. 
Then the conductors of the coordinate monoids of $S$ are bounded by $k$.
In particular for every $j\in \{1,\dots, n-1\}$ we have 
\begin{align*}
\mtt{n}(S_j) \leq  k-\mtt{g}_j.
\end{align*}
Thus, it follows from Lemma~\ref{L:sporadicelements} that
\begin{align}\label{A:toberearranged}
(k-\mtt{g}_1)\cdots (k-\mtt{g}_{n-1}) \leq k^{n-1} - \mtt{g}_M(S).
\end{align}
After expanding the left hand side of (\ref{A:toberearranged}) by using the elementary symmetric polynomials, and then reorganizing the result, we find our desired inequality in (\ref{A:genera}). 
This finishes the proof of our proposition.
\end{proof}

The proof of the following corollary readily follows from Proposition~\ref{P:generacomparison}.

\begin{Corollary}
Let $S$ be a unipotent numerical submonoid of $\mbf{P}(3,\N)$.
Let $\mtt{g}_1,\mtt{g}_2$ denote the genera of the coordinate submonoids $S_1,S_2$, respectively.
Then we have 
\begin{align*}
\mtt{g}_1+\mtt{g}_2\leq \mtt{g}_{\mbf{P}(3,\N)}(S) \leq  k(\mtt{g}_1+\mtt{g}_2)-\mtt{g}_1\mtt{g}_2,
\end{align*}
where $k$ is the generating number of $S$ relative to $\mbf{P}(3,\N)$.
\end{Corollary}

\section{The Thick and the Thin Unipotent Numerical Groups}\label{S:Thick}

In this section, we follow the notational conventions of the previous section.
In particular, we will not distinguish between $\mbf{P}(n,\N)$ and $\N^{n-1}$.
The generalized monoids that attain the lower bound in Proposition~\ref{P:generacomparison} are rather special. Indeed, if the equality $\mtt{g}_1+\cdots +\mtt{g}_{n-1} = \mtt{g}_M(S)$ holds, then the elements of $M\setminus S$ are 
all contained in the union $\bigcup_{j=1}^{n-1} S_j$. 
This means that the set $M\setminus \bigcup_{j=1}^{n-1} M_j$ is a subsemigroup of $S$.

\begin{Definition}
Let $S$ be a unipotent numerical monoid in $M$ with coordinate submonoids $S_1,\dots, S_{n-1}$. 
Let $\mtt{g}_1,\dots, \mtt{g}_{n-1}$ denote the genera of $S_1,\dots,S_{n-1}$, respectively. 
Let $\mtt{n}_1,\dots, \mtt{n}_{n-1}$ denote the sporadicities of $S_1,\dots, S_{n-1}$, respectively. 
If $\sum_{j=1}^{n-1} \mtt{g}_{j} = \mtt{g}_M(S)$ holds, then we call $S$ a {\em thick (unipotent numerical) monoid in $M$}. 
If $\prod_{j=1}^{n-1} \mtt{n}_j= \mtt{n}_M(S)$ holds, then we call $S$ a {\em thin (unipotent numerical) monoid in $M$}.
\end{Definition}

Clearly, if $n=2$, then every unipotent numerical monoid is simultaneously a thin monoid and a thick monoid.
In our next result, we prove that, for $n\geq 3$, the intersection of the families of thin and thick monoids is trivial.
\begin{Theorem}
For $n\geq 3$, a unipotent numerical monoid $S$ in $M$ is simultaneously a thin monoid and a thick monoid if and only if $S= M$.
\end{Theorem}

\begin{proof}

($\Leftarrow$) If $S=M$, then we have $\mtt{g}_M(S)=0$ and $\mtt{n}_M(S)=1$.
Hence, the conditions $\sum_{j=1}^{n-1} \mtt{g}_{j} = \mtt{g}_M(S)$ and $\prod_{j=1}^{n-1} \mtt{n}_j= \mtt{n}_M(S)$
are automatically satisfied. In other words, if $S=M$, then $S$ is both a thin and a thick monoid. 
\medskip

($\Rightarrow$)
The following simple lemma will be useful.
\begin{Lemma}\label{L1:thinthick}
Let $k\geq 2$. Let $\lambda_1\leq \dots \leq \lambda_m$ be a weakly increasing sequence of nonnegative integers such
$\lambda_m \leq k-1$. 
Then we have $\frac{1}{k}\left(\sum_{i=1}^m \lambda_i \right) \leq \sum_{i=2}^m \lambda_i$. 
\end{Lemma}
\begin{proof}
We will show that $\sum_{i=1}^m \lambda_i  \leq k \left( \sum_{i=2}^m \lambda_i \right)$.
Let us look closely at the right hand side of the following decomposition:
\begin{align*}
k\left(\sum_{i=2}^m \lambda_i\right) = (k-1) \left( \sum_{i=2}^m \lambda_i \right)  + \left( \sum_{i=2}^m \lambda_i \right).
\end{align*}
Since $k\geq 2$, we have  $(k-1) \left( \sum_{i=2}^m \lambda_i \right) \geq   \sum_{i=2}^m \lambda_i$.
Since $ \sum_{i=2}^m \lambda_i \geq \lambda_1$, we have 
\begin{align*}
(k-1) \left( \sum_{i=2}^m \lambda_i \right)  + \left( \sum_{i=2}^m \lambda_i \right) \geq 
\left(\sum_{i=2}^m \lambda_i \right) + \lambda_1 = \sum_{i=1}^m \lambda_i.
\end{align*}
This finishes the proof of our lemma.
\end{proof}

Now let $S\subset M = \mbf{P}(n,\N)$ be a unipotent numerical semigroup which is both thin and thick at the same time. 
Let $k$ denote generating number of $S$. 
Our goal is show that $k=1$. 
Towards a contradiction, we proceed with the assumption that $k\geq 2$.

Since $S$ is a thick monoid, it follows from the general identity 
\[
\mtt{n}_M(S) + \mtt{g}_M(S) = \mtt{c}_M(S) = k^{n-1}
\] 
that the sporadicity of $S$ relative to $M$ is given by 
$\mtt{n}_M(S) = k^{n-1} - \sum_{i=1}^{n-1} \mtt{g}_i$, where $k$ denotes the conductor of $S$ relative to $M$.
Since $\mtt{g}_i + \mtt{n}_i \leq k$ for every $i\in \{1,\dots, n-1\}$, and since $S$ is a thin monoid, we obtain the following inequality, 
\begin{align}\label{A:thinthick}
k^{n-1} - \sum_{i=1}^{n-1} \mtt{g}_i = \mtt{n}_M(S) = \prod_{i=1}^{n-1}\mtt{n}_i \leq \prod_{i=1}^{n-1} (k -\mtt{g}_i).
\end{align}

The following lemma will provide us with the last step of our proof.
\begin{Lemma}\label{L2:thinthick}
Let $m,k\geq 2$. Let $\lambda_1\leq \dots \leq \lambda_m$ be a weakly increasing sequence of nonnegative integers such
$\lambda_m \leq k-1$. 
If the inequality 
\begin{align}\label{A:keyinequality}
k^{m} - \sum_{i=1}^{m} \lambda_i \leq \prod_{i=1}^{m} (k -\lambda_i)
\end{align} 
is satisfied, then we have 
$\lambda_1 = \cdots = \lambda_m=0$. 
\end{Lemma}
\begin{proof}
We prove this direction by using induction on $m$. 
The base case is when $m=2$.
Then our inequality is given by $k^2 - (\lambda_1+\lambda_2) \leq (k-\lambda_1)(k-\lambda_2)$.
After reorganizing it, we find that $(\lambda_1+\lambda_2)(k-1) \leq \lambda_1\lambda_2$. 
Since $k\geq 2$ and $\lambda_i \leq k-1$ for $i\in \{1,2\}$, the last inequality holds if and only if $\lambda_1=\lambda_2=0$.
We now assume that our claim holds true for every $l\in \{2,\dots, m-1\}$, and we proceed to prove it for $l=m$.
To this end, we divide both sides of (\ref{A:keyinequality}) by $k$, and reorganize it: 
$
k^{m-1} \leq \frac{1}{k} \prod_{i=1}^{m} (k -\lambda_i) +  \frac{1}{k}\sum_{i=1}^{m} \lambda_i.
$
Note that $(k-\lambda_1)/k \leq 1$. 
Then by applying Lemma~\ref{L1:thinthick}, we obtain 
\begin{align*}
k^{m-1} \leq \prod_{i=2}^{m} (k -\lambda_i) +   \sum_{i=2}^m \lambda_i.
\end{align*} 
After relabeling our new inequality and then applying the induction assumption, we see that $\lambda_2 =\cdots = \lambda_m=0$.
But $\lambda_1$ satisfies $0\leq \lambda_1\leq \lambda_2$. 
Hence, the proof of our lemma is finished. 
\end{proof}
We are now ready to finish the proof of our theorem. 
We apply Lemma~\ref{L2:thinthick} to the inequality (\ref{A:thinthick}) with $\lambda_i = \mtt{g}_i$ for $i\in \{1,\dots, n-1\}$.
Then we see that $\mtt{g}_1=\cdots = \mtt{g}_{n-1}=0$. 
But if the genera of all coordinate submonoids of $S$ are 0, then the genus of $S$ relative to $M$ must be 0. 
This means that $S=M$. 
Hence, the proof of our theorem is finished. 
\end{proof}

We proceed to show that the Unipotent Wilf Conjecture holds for the thick monoids when $n$ is at least 3.

\begin{Theorem}\label{T:UWCforthick}
Let $n\geq 3$, and let $M:=\mbf{P}(n,\N)$.
If $S$ is a thick unipotent numerical monoid in $S$, 
then the Unipotent Wilf Conjecture~\ref{C:secondconjecture} holds true for $S$. 
\end{Theorem}

\begin{proof}
Let $S_i$ ($i\in \{1,\dots, n-1\}$) denote the $i$-th coordinate monoid of $S$. 
Let $\mtt{e}_i$ and $\mtt{n}_i$ denote the embedding dimension and the sporadicity of $S_i$, respectively. 
Clearly, every generator of $S_i$ is a generator of $S$. 
Hence, we have the basic estimate, $\sum_{i=1}^{n-1} \mtt{e}_i \leq \mtt{e}(S)$.
Let $k$ denote the generating number of $S$ relative to $M$. 
Without loss of generality, we assume that $k\geq 2$. 
Hence, one of the numbers $\mtt{e}_1,\dots, \mtt{e}_{n-1}$ must be at least two. 
It follows that 
\begin{align} \label{A:basicestimatefore}
n\ \leq \ \sum_{i=1}^{n-1} \mtt{e}_i \leq \mtt{e}(S).
\end{align}
At the same time, since $S$ is a thick monoid, the number of sporadic elements is found by the following formula: 
\begin{align} \label{A:basicestimateforn}
\mtt{n}_M(S) = k^{n-1} - \sum_{i=1}^{n-1} \mtt{g}_i.
\end{align}
Here, $k^{n-1}$ is the cardinality of the enveloping hypercube of $M\setminus S$. 
By putting (\ref{A:basicestimatefore}) and (\ref{A:basicestimateforn}) together, we obtain $(k^{n-1} - \sum_{i=1}^{n-1} \mtt{g}_i) n$ as a lower bound for $\mtt{e}(S) \mtt{n}_M(S)$.
Notice that 
\begin{align*}
(k^{n-1} - \sum_{i=1}^{n-1} \mtt{g}_i) n = (n-1)k^{n-1} +k^{n-1} - n \left( \sum_{i=1}^{n-1} \mtt{g}_i\right).
\end{align*}
This means that if the inequality $0\leq k^{n-1} - n \left( \sum_{i=1}^{n-1} \mtt{g}_i\right)$ holds, then the proof of our claim will follow. 
Since the conductor of $S_i$ ($i\in \{1,\dots, n-1\}$) is at most $k$, we see that $\mtt{g}_i \leq k-1$. 
Therefore, we see that 
\begin{align*}
k^{n-1} - n(n-1)(k-1) \leq k^{n-1} - n \left( \sum_{i=1}^{n-1} \mtt{g}_i\right).
\end{align*}
To analyze the lower bound, we view it as a function on $\Z^2$,  
\begin{align*}
f : \Z\times \Z &\longrightarrow \Z \\
(n,k) &\longmapsto  k^{n-1} - n(n-1)(k-1).
\end{align*}
We are interested in the values $f(n,k)$, where $k\geq 2$ and $n\geq 3$. 
For such a pair $(n,k)$, it is easy to see by using Calculus that 
\begin{align*}
f(n,k) < 0 \qquad\text{if}\qquad  (n,k)\in \{ (3,2),(3,3),(3,4),(4,2),(5,2)\}.
\end{align*}
But it is easy to see by using Calculus that, for every $k\geq 5$ and $n\geq 3$, the function $(k,n)\mapsto k^{n-1} - n(n-1)(k-1)$ takes values in $\Z_+$. 
Hence, for such $(n,k)$, our claim holds true. 
For $ (n,k)\in \{(3,2),(3,3),(3,4),(4,2),(5,2)\}$, it is easy to check (by hand) that the Unipotent Wilf Conjecture holds true. 
This finishes the proof of our theorem.
\end{proof}

\begin{Example}
In this example, for the convenience of the reader, we will check on an example, where $(n,k)=(5,2)$, that the Unipotent Wilf Conjecture holds. 
Recall that we identify $\mbf{P}(5,\N)$ with $\N^4$. 
Let $S$ be the thick monoid in $\N^4$ defined by 
\begin{align*}
S =\langle (1,0,0,0),(0,0,0,1),(1,1,1,1), A:\ A \in \mc{P}_2\rangle. 
\end{align*}
Then it is easy to verify that the minimal generating set is given by 
\[
\left\{
\begin{array}{lll}
(1,0,0,0), & (0,0,0,1), & (1,1,1,1),\\
(1,0,0,1), & (0,2,0,0), & (0,3,0,0),\\  
(0,1,1,2), & (0,0,1,2), & (0,1,0,2), \\ 
(0,2,1,0), & (0,3,1,0), & (0,0,2,0), \\
(0,0,3,0), & (0,1,2,0), & (0,1,3,0), \\
(2,1,1,0), & (2,0,1,0), & (2,1,0,0)
\end{array} \right\}.
\]
In particular, we see that $\mtt{e}(S) = 18$. 
The elements of $S$ that are contained in the enveloping hypercube of $M\setminus S$ are 
\begin{align*}
(0,0,0,0),(1,0,0,0),(0,0,0,1),(1,0,0,1),(1,1,1,1).
\end{align*}
Hence, the sporadicity of $S$ relative to $M$ is $5$.
It follows that $\mtt{e}(S) \cdot \mtt{n}_M(S) = 90$ while $4\mtt{c}_M(S)=64$ since the generating number is 2.
\end{Example}

\medskip

Next, we will consider the thin monoids. 
First we introduce the well-known notion of a multiplicity of a numerical monoid.
If $T$ is a numerical monoid, then the number
\begin{align*}
\min\{ a: a\in T\setminus \{0\}\}
\end{align*}
is called the {\em multiplicity} of $T$.

\begin{Lemma}\label{L:embeddingdimofthin}
Let $S$ be a thin unipotent numerical monoid in $M:=\mbf{P}(3,\N)$ with coordinate monoids $S_1$ and $S_2$. 
Then the embedding dimension of $S$ is given by
\begin{align*}
\mtt{e}(S):= \mtt{e}_1 +\mtt{e}_2 + (\mtt{m}_1-1)\mtt{m}_2 + \mtt{m}_1(\mtt{m}_2-1), 
\end{align*}
where $\mtt{m}_i$ ($i\in \{1,2\}$) is the multiplicity of $S_i$. 
\end{Lemma}

\begin{proof}
Let $A$ denote the enveloping square of $M\setminus S$.
The minimal generators of a unipotent numerical monoid $S$ in $M$ are always contained in the union of 
the enveloping hypercube of $M\setminus A$ and the set of minimal generators fo the $k$-th fundamental monoid of $S$,
where $k$ is the conductor of $S$ relative to $M$. 
Since $S$ is a thin monoid, the set $A\setminus S_1\cup S_2$ does not share any elements with $\mc{A}$. 
In other words, all minimal generators of $S$ that are contained in $A$ are actually contained in the coordinate submonoids. 
In fact, more is true. 
The minimal generators of the coordinate monoids are always among the minimal generators for $S$. 
Let us denote the set of minimal generators of $S$ that are contained in $S_1\cup S_2$ by $\mc{A}$. 
Hence, we have $|\mc{A}| =\mtt{e}_1+\mtt{e}_2$.
Thus, it remains to characterize those minimal generators of $S$ that are contained in $\mc{P}_k \setminus S_1\cup S_2$,
where $\mc{P}_k$ is the minimal generating set for $P(k)$.

Let $(a_1,\dots, a_{\mtt{m}_1}) := (1,2,\dots, \mtt{m}_1)$ and $(b_1,\dots, b_{\mtt{m}_2}) := (1,2,\dots, \mtt{m}_2)$.
It is not difficult to see that each pair of the form $(k+s,a_i)$, where $i\in \{1,2,\dots, \mtt{m}_1-1\}$ and $s\in \{0,1,\dots, \mtt{m}_2-1\}$ 
is one of the minimal generators of $S$ contained in $\mc{P}_k \setminus S_1\cup S_2$.
Let $\mc{A}_1$ denote the set of all such pairs. 
Clearly, the cardinality of $\mc{A}_1$ is given by $(\mtt{m}_1-1)\mtt{m}_2$.
Likewise, it is not difficult to see that each pair of the form $(b_j,k+s)$, where $j\in \{1,2,\dots, \mtt{m}_2-1\}$ and $s\in \{0,1,\dots, \mtt{m}_1-1\}$ is a minimal generator in $\mc{P}_k \setminus S_1\cup S_2$ as well. 
Let $\mc{A}_2$ denote the set of all such pairs. 
The cardinality of $\mc{A}_2$ is given by $\mtt{m}_1(\mtt{m}_2-1)$.
Then, the union $\mc{A}_1\sqcup \mc{A}_2$ is a subset of the minimal generating set of $S$. 
We now observe that, if we translate $\mc{A}_1$ (resp. $\mc{A}_2$) to the origin by subtracting $(k,0)$ from its elements, 
then the resulting set becomes a generating set for $\N^2$. 
Likewise the translated set $\mc{A}_2 - (0,k)$ generates $\N^2$. 
But this means that the union $\mc{A}\sqcup \mc{A}_1\sqcup \mc{A}_2 $ generates $S$. 
Since this union is already a subset of the minimal generating set of $S$, it must be equal to it. 
Finally, we notice that $\mtt{e}_1 +\mtt{e}_2 +  (\mtt{m}_1-1)\mtt{m}_2 + \mtt{m}_1( \mtt{m}_2-1) = |\mc{A}\sqcup \mc{A}_1\sqcup \mc{A}_2|$. This finishes the proof of our proposition.
\end{proof}

We will illustrate the proof of our previous proposition on an example.

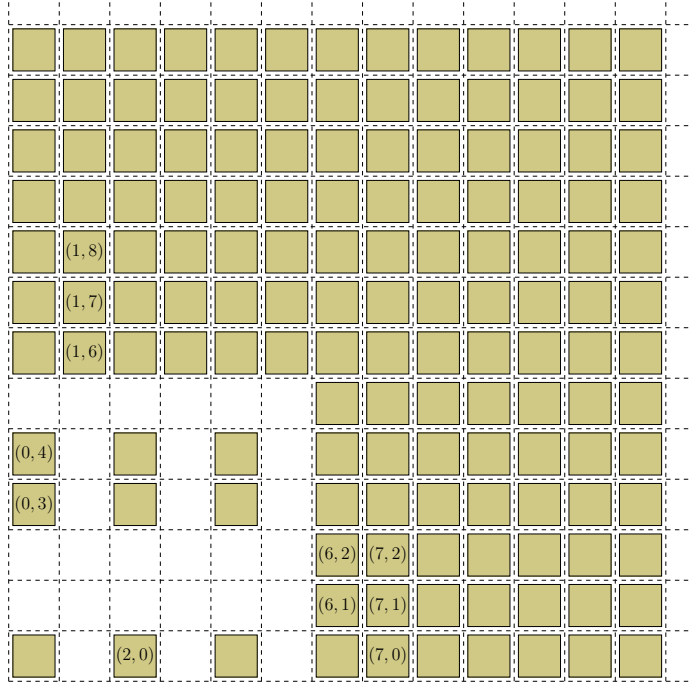
\begin{figure}[htp]
\begin{center}
\scalebox{.56}{
\begin{tikzpicture}[scale=1.2]

\foreach \i in {0,2,4,6,7,8,9,10,11,12} {\draw (\i,0)  node[fill=olive!45!, minimum size=1cm,draw] {};};
\foreach \i in {0,2,4,6,7,8,9,10,11,12} {\draw (\i,3)  node[fill=olive!45!, minimum size=1cm,draw] {};};
\foreach \i in {0,2,4,6,7,8,9,10,11,12} {\draw (\i,4)  node[fill=olive!45!, minimum size=1cm,draw] {};};

\foreach \i in {0,1,...,12} {\draw (\i,6)  node[fill=olive!45!, minimum size=1cm,draw] {};};
\foreach \i in {0,1,...,12} {\draw (\i,7)  node[fill=olive!45!, minimum size=1cm,draw] {};};
\foreach \i in {0,1,...,12} {\draw (\i,8)  node[fill=olive!45!, minimum size=1cm,draw] {};};
\foreach \i in {0,1,...,12} {\draw (\i,9)  node[fill=olive!45!, minimum size=1cm,draw] {};};
\foreach \i in {0,1,...,12} {\draw (\i,10)  node[fill=olive!45!, minimum size=1cm,draw] {};};
\foreach \i in {0,1,...,12} {\draw (\i,11)  node[fill=olive!45!, minimum size=1cm,draw] {};};
\foreach \i in {0,1,...,12} {\draw (\i,12)  node[fill=olive!45!, minimum size=1cm,draw] {};};
\foreach \i in {6,...,12} {\draw (\i,1)  node[fill=olive!45!, minimum size=1cm,draw] {};};
\foreach \i in {6,...,12} {\draw (\i,2)  node[fill=olive!45!, minimum size=1cm,draw] {};};
\foreach \i in {6,...,12} {\draw (\i,5)  node[fill=olive!45!, minimum size=1cm,draw] {};};

\node at (2,0) {$(2,0)$};
\node at (0,3) {$(0,3)$};
\node at (0,4) {$(0,4)$};
\node at (7,0) {$(7,0)$};
\node at (1,6) {$(1,6)$};
\node at (1,7) {$(1,7)$};
\node at (1,8) {$(1,8)$};
\node at (6,1) {$(6,1)$};
\node at (7,1) {$(7,1)$};
\node at (6,2) {$(6,2)$};
\node at (7,2) {$(7,2)$};
\foreach \i in {0,...,13} {\draw [dashed] (\i-.5,-.5) -- (\i-.5,13);};
\foreach \i in {0,...,13} {\draw [dashed] (-.5,\i-0.5) -- (13,\i-0.5);};

\end{tikzpicture}
}
\end{center}
\caption{A thin monoid.}
\label{F:thin}
\end{figure}

\begin{Example}
Let $S$ denote the unipotent numerical monoid in $\mbf{P}(3,\N)$ that is depicted in Figure~\ref{F:thin}.
As usual, the shaded boxes in this figure represent the elements of $S$. 
Let $S_1$ and $S_2$ denote the coordinate monoids of $S$. 
Let $\mc{A} := \{ (0,3), (0,4), (2,0),(7,0)\}$. 
Then $\mc{A}$ consists of the minimal generators of $S$ contained in $S_1\cup S_2$. 
The following sets give the remaining part of the minimal generating set of $S$:
\begin{align*}
\mc{A}_1 = \{ (1,6), (1,7), (1,8)\}\qquad\text{and}\qquad
\mc{A}_2 = \{ (6,1), (6,2),(7,1),(7,2)\}.
\end{align*}
It is easy to check that $\mc{A}_1\sqcup \mc{A}_2 \sqcup \mc{A}$ is a partitioning of the minimal generating set of $S$. 
\end{Example}

Next, we are now going to extend the conclusion of Lemma~\ref{L:embeddingdimofthin} to higher dimensions.

\begin{Lemma}\label{L:embeddingdimofthin2}
Let $S$ be a thin unipotent numerical monoid in $M:=\mbf{P}(n,\N)$ with coordinate monoids $S_1,\dots, S_{n-1}$. 
We have the following identity for the embedding dimension of $S$:
\begin{align*}
\mtt{e}(S):= \sum_{i=1}^{n-1} \mtt{e}_i +  \sum_{\substack{I\subseteq \{1,\dots, n-1\},\\ |I|\geq 2}} 
\left( \sum_{i\in I} \mtt{m}_i \prod_{j\in I\setminus\{ i\} } (\mtt{m}_j-1) \right).
\end{align*}
\end{Lemma}

\begin{proof}
The proof follows the line of reasoning that we used in Lemma~\ref{L:embeddingdimofthin}.
This time we have $n-1$ stages of computation instead of two. 
Since the details of the specific computations are very similar to the details of the proof of Lemma~\ref{L:embeddingdimofthin},
we provide the final formulas only. 
For a nonempty subset $I \subseteq \{1,\dots, n-1\}$, let $\N^{n-1}_I$ denote the coordinate subspace consisting of 
vectors $(v_1,\dots, v_{n-1})\in \N^n$ such that $v_i=0$ for every $i\notin I$. 

In the first stage of our calculation, we count the minimal generators that are contained in the sets 
$\N^{n-1}_I \cap S$, where $I$ varies over all 1-element subsets of $\{1,\dots, n-1\}$.
In other words, we count the number of minimal generators in the union, $S_1\cup \cdots \cup S_{n-1}$. 
This is given by the sum, 
\begin{align}\label{A:firstsum}
\sum_{i=1}^{n-1} \mtt{e}_i.
\end{align}

In the second step, we calculate the number of minimal generators that are contained in the sets  
$\N^{n-1}_I \cap S$, where $I$ varies over all 2-element subsets of $\{1,\dots, n-1\}$.
This number is given by 
\begin{align}\label{A:secondsum}
\sum_{\{i,j\}\subset \{1,\dots, n-1\}} (\mtt{m}_i-1) \mtt{m}_j.
\end{align} 

More generally, in the $r$-th stage, where $r\in \{3,\dots, n-1\}$, 
we calculate the number of minimal generators that are contained in the sets  
$\N^{n-1}_I \cap S$, where $I$ varies over all $r$-element subsets of $\{1,\dots, n-1\}$.
This number is given by 
\begin{align}\label{A:rthsum}
\sum_{\substack{I\subset \{1,\dots, n-1\} \\ |I|=r}} \sum_{i\in I} \mtt{m}_i \prod_{\substack{j \in I\setminus \{ i\}}} (\mtt{m}_j-1).
\end{align} 
Putting (\ref{A:firstsum})--(\ref{A:rthsum}) together we obtain our formula for the number of minimal generators for $S$. 
\end{proof}

\begin{Theorem}\label{T:realformula}
Let $S$ be a thin unipotent numerical monoid in $M:=\mbf{P}(n,\N)$ with coordinate monoids $S_1,\dots, S_{n-1}$. 
Then we have the following lower bound for the embedding dimension of $S$:
\begin{align*}
(n-1) \prod_{i=1}^{n-1} \mtt{e}_i\  \leq \ \mtt{e}(S).
\end{align*}
\end{Theorem}

\begin{proof}
Since the embedding dimension of a numerical monoid is less than or equal to its multiplicity, 
by Lemma~\ref{L:embeddingdimofthin2}, we have 
\begin{align}\label{A:someinequality}
 \sum_{i=1}^{n-1} \mtt{e}_i +  \sum_{\substack{I\subseteq \{1,\dots, n-1\},\\ |I|\geq 2}} 
\left( \sum_{i\in I} \mtt{e}_i \prod_{j\in I\setminus\{ i\} } (\mtt{e}_j-1) \right) \ \leq \ \mtt{e}(S).
\end{align}
We will use induction on $n$ to prove that the right hand side of (\ref{A:someinequality}) is given by $(n-1) \prod_{i=1}^{n-1} \mtt{e}_i$.
For $n=2$, we already proved our claim in Lemma~\ref{L:embeddingdimofthin}.
We assume that our claim holds for $n-1$.
Let us proceed to prove it for $n$. 
We split the sum 
$\sum_{\substack{I\subseteq \{1,\dots, n\},\\ |I|\geq 2}} \left( \sum_{i\in I} \mtt{e}_i \prod_{j\in I\setminus\{ i\} } (\mtt{e}_j-1) \right)$
 into three parts, denoting them by $A$ and $B$:
\begin{align*}
\underbrace{\sum_{\substack{I\subseteq \{1,\dots, n-1\},\\ |I| \geq 2}} 
\left( \sum_{i\in I} \mtt{e}_i \prod_{j\in I\setminus\{ i\} } (\mtt{e}_j-1) \right)}_{A} + 
\underbrace{\sum_{\substack{I\subseteq \{1,\dots, n\},\\ n\in I,\ |I| \geq 2}} 
\left( \sum_{i\in I} \mtt{e}_i \prod_{j\in I\setminus\{ i\} } (\mtt{e}_j-1) \right)}_{B}.   
\end{align*}
Then we have 
\begin{align*}
\sum_{i=1}^n \mtt{e}_i + A + B &= 
\left(\sum_{i=1}^{n-1} \mtt{e}_i + A\right) + \mtt{e}_n + B  \\
&= (n-1) \prod_{i=1}^{n-1} \mtt{e}_i + \mtt{e}_n + B \qquad (\text{by induction}).
\end{align*}
Let $I$ be an indexing subset in $B$. 
Then, denoting them by $B_1(I)$ and $B_2(I)$, we split the inner sum of $B$ corresponding to $I$ into two parts:
\begin{align*}
\sum_{i\in I} \mtt{e}_i \prod_{j\in I\setminus\{ i\} } (\mtt{e}_j-1)  =
\underbrace{\mtt{e}_n \prod_{j\in I\setminus\{n \} } (\mtt{e}_j-1)}_{B_1(I)} +
\underbrace{(\mtt{e}_n-1)\left( \sum_{i\in I\setminus \{n\}} \mtt{e}_i \prod_{j\in I\setminus\{ i,n\} } (\mtt{e}_j-1) \right)}_{B_2(I)}.
\end{align*}
Note that $I\setminus \{n\}$ may have one element only. 
We distinguish these possibilities in the sums $\sum_{\substack{I\subseteq \{1,\dots, n\} \\ n\in I}} B_2(I)$:
\begin{align*}
\sum_{\substack{I\subseteq \{1,\dots, n\} \\ n\in I,\ |I|\geq 2}} B_2(I) 
&= (\mtt{e}_n-1)\left( \sum_{i=1}^{n-1} \mtt{e}_i 
+ \sum_{\substack{I'\subseteq \{1,\dots, n-1\} \\ |I'|\geq 2}}\sum_{i\in I'} \mtt{e}_i \prod_{j\in I'\setminus\{ i\} } (\mtt{e}_j-1) \right).
\end{align*}
Now our induction hypothesis is applicable to the sum in the parenthesis.
We obtain, 
\begin{align*}
\sum_{\substack{I\subseteq \{1,\dots, n\}, \\ n\in I,\ |I|\geq 2}} B_2(I) &= 
(\mtt{e}_n-1) (n-1) \prod_{i=1}^{n-1} \mtt{e}_i. 
\end{align*}
In summary, our sum $(n-1) \prod_{i=1}^{n-1} \mtt{e}_i + \mtt{e}_n + B$ becomes 
\begin{align*}
(n-1) \prod_{i=1}^{n-1} \mtt{e}_i + \mtt{e}_n + \sum_{\substack{I\subseteq \{1,\dots, n\},\\ n\in I,\ |I|\geq 2}} B_1(I) + 
(\mtt{e}_n-1) (n-1) \prod_{i=1}^{n-1} \mtt{e}_i.
\end{align*}
We perform the obvious cancellation. 
Then we carefully expand the summation involving $B_1(I)$:
\begin{align}\label{A:nowuse}
 (n-1) \prod_{i=1}^{n} \mtt{e}_i  +  \mtt{e}_n +
\mtt{e}_n\sum_{\substack{I\subseteq \{1,\dots, n-1\},\\ |I|\geq 1}}  \prod_{j\in I\ } (\mtt{e}_j-1).
\end{align}
But it is a simple combinatorial exercise to show that 
\begin{align*}
\prod_{i=1}^{n-1} (1+x_i) = 1 + \sum_{ \substack{I\subseteq \{1,\dots, n-1\},\\ I\neq \emptyset}} \prod_{i\in I} x_i.
\end{align*}
In particular, if we set $x_i =\mtt{e}_i-1$ in (\ref{A:nowuse}), then we see that our sum becomes 
\begin{align}\label{A:nowuse2}
 (n-1) \prod_{i=1}^{n} \mtt{e}_i  +  \mtt{e}_n + \mtt{e}_n\left( \prod_{i=1}^{n-1} \mtt{e}_i-1\right).
\end{align}
A final round of cancellations show that (\ref{A:nowuse2}) is equal to $n\prod_{i=1}^n \mtt{e}_i$. 
This finishes the proof of our theorem.
\end{proof}

In the introduction we mentioned that for certain unipotent numerical monoids, 
the ingredients of the Generalized Wilf Conjecture are the same as those of the Unipotent Wilf Conjecture. 
In our next result we have one of these situations.

\begin{Theorem}\label{T:SUWCforthin}
We follow the notation of our previous theorem. 
Let $S$ be a thin monoid in $M:=\mbf{P}(n,\N)$. 
Let $k$ denote the generating number of $S$.
We assume that the ordinary Wilf conjecture holds for the coordinate monoids $S_1,\dots, S_{n-1}$. 
If the identity $\prod_{i=1}^{n-1} \mtt{c}_i=k^{n-1}$ holds for the conductors of the $S_i$'s, then the Unipotent Wilf Conjecture holds for $S$. 
Furthermore, in this case, the Unipotent Wilf Conjecture is equivalent to the Generalized Wilf Conjecture.
\end{Theorem}

Before we start our proof, let us point out that our conclusion of the equivalence of the Unipotent Wilf Conjecture to the Generalized Wilf Conjecture 
is independent of the ordinary Wilf conjecture.

\begin{proof}
Since we know that $\mtt{c}_i$ ($i\in \{1\dots, n-1\}$) is at most $k$, 
the identity $\prod_{i=1}^{n-1} \mtt{c}_i=k^{n-1}$ holds if and only if the equalities $\mtt{c}_1=\cdots = \mtt{c}_{n-1} = k$ hold. 
This means that the hole set of $S$ in $\N^{n-1}$ contains the vector $(k-1,\dots, k-1)$. 
It follows that $c(S)=\mtt{c}_M(S)$ and $n(S) = \mtt{n}_M(S)$. 
In other words, the Unipotent Wilf Conjecture is equivalent to the Generalized Wilf Conjecture under our assumptions. 

We are now ready to prove our main claim.
By Theorem~\ref{T:realformula}, we have 
$(n-1) \prod_{i=1}^{n-1}\mtt{e}_i \leq \mtt{e}_M(S)$.
Multiplying both sides of this inequality with $\mtt{n}_M(S)$,
and using the ordinary Wilf conjecture, yield the following inequalities: 
\begin{align*}
(n-1) \ \prod_{i=1}^{n-1} \mtt{c}_i = (n-1)k^{n-1} \ \leq \ (n-1) \ \prod_{i=1}^{n-1}\mtt{e}_i \mtt{n}_i  \ \leq \ \mtt{e}_M(S) \mtt{n}_M(S).
\end{align*}
This finishes the proof of our theorem.
\end{proof}

We close our paper by an example of a thin monoid for which the 
Generalized Wilf Conjecture and the Unipotent Wilf Conjecture assert different inequalities.

\begin{Example}
Let $S$ denote the thin monoid depicted in Figure~\ref{F:thin2}, 
where the shaded boxes represent the elements of $S$. 
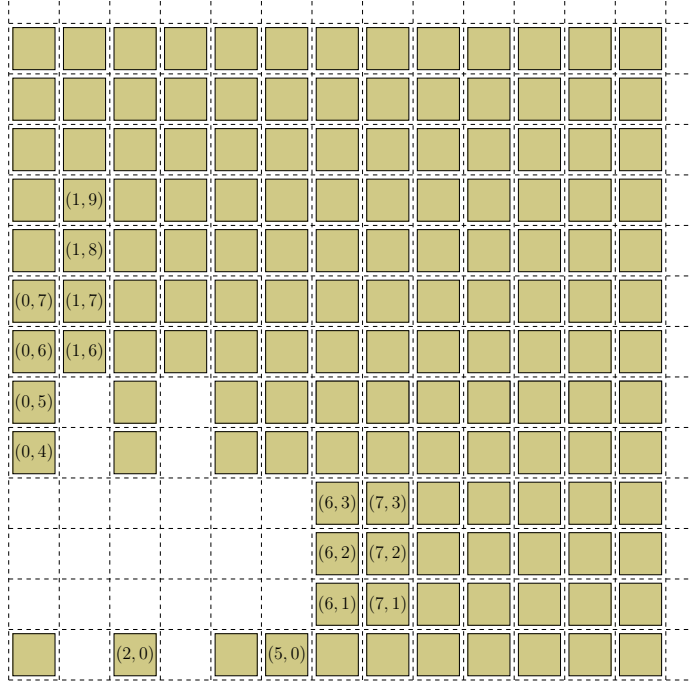
\begin{figure}[htp]
\begin{center}
\scalebox{.56}{
\begin{tikzpicture}[scale=1.2]

\foreach \i in {0,2,4,5,6,7,8,9,10,11,12} {\draw (\i,0)  node[fill=olive!45!, minimum size=1cm,draw] {};};
\foreach \i in {0,2,4,5,6,7,8,9,10,11,12} {\draw (\i,4)  node[fill=olive!45!, minimum size=1cm,draw] {};};
\foreach \i in {0,2,4,5,6,7,8,9,10,11,12} {\draw (\i,5)  node[fill=olive!45!, minimum size=1cm,draw] {};};
\foreach \i in {6,7,8,9,10,11,12} {\draw (\i,3)  node[fill=olive!45!, minimum size=1cm,draw] {};};
\foreach \i in {0,1,...,12} {\draw (\i,6)  node[fill=olive!45!, minimum size=1cm,draw] {};};
\foreach \i in {0,1,...,12} {\draw (\i,7)  node[fill=olive!45!, minimum size=1cm,draw] {};};
\foreach \i in {0,1,...,12} {\draw (\i,8)  node[fill=olive!45!, minimum size=1cm,draw] {};};
\foreach \i in {0,1,...,12} {\draw (\i,9)  node[fill=olive!45!, minimum size=1cm,draw] {};};
\foreach \i in {0,1,...,12} {\draw (\i,10)  node[fill=olive!45!, minimum size=1cm,draw] {};};
\foreach \i in {0,1,...,12} {\draw (\i,11)  node[fill=olive!45!, minimum size=1cm,draw] {};};
\foreach \i in {0,1,...,12} {\draw (\i,12)  node[fill=olive!45!, minimum size=1cm,draw] {};};
\foreach \i in {6,...,12} {\draw (\i,1)  node[fill=olive!45!, minimum size=1cm,draw] {};};
\foreach \i in {6,...,12} {\draw (\i,2)  node[fill=olive!45!, minimum size=1cm,draw] {};};
\foreach \i in {6,...,12} {\draw (\i,5)  node[fill=olive!45!, minimum size=1cm,draw] {};};

\node at (2,0) {$(2,0)$};
\node at (0,4) {$(0,4)$};
\node at (0,5) {$(0,5)$};
\node at (0,6) {$(0,6)$};
\node at (0,7) {$(0,7)$};
\node at (5,0) {$(5,0)$};
\node at (1,6) {$(1,6)$};
\node at (1,7) {$(1,7)$};
\node at (1,8) {$(1,8)$};
\node at (1,9) {$(1,9)$};
\node at (6,1) {$(6,1)$};
\node at (7,1) {$(7,1)$};
\node at (6,2) {$(6,2)$};
\node at (6,3) {$(6,3)$};
\node at (7,2) {$(7,2)$};
\node at (7,3) {$(7,3)$};
\foreach \i in {0,...,13} {\draw [dashed] (\i-.5,-.5) -- (\i-.5,13);};
\foreach \i in {0,...,13} {\draw [dashed] (-.5,\i-0.5) -- (13,\i-0.5);};

\end{tikzpicture}
}
\end{center}
\caption{A thin monoid that admits two different Wilf conjectures.}
\label{F:thin2}
\end{figure}

The boxes with coordinates in them correspond to the minimal generators of $S$.
Then, we have $\mtt{e}(S)= 16$, $\mtt{n}_M(S) =12$, and $\mtt{r}_M(S)=6$. 
The Unipotent Wilf Conjecture predicts correctly the following inequality:
\begin{align*}
2\cdot \mtt{c}_M(S)=72\ \leq \ 16\cdot 12 = \mtt{e}(S)\cdot \mtt{n}_M(S).
\end{align*}
 At the same time, the Generalized Wilf Conjecture predicts correctly the following inequality:
\begin{align*}
2\cdot c(S)= 68\ \leq \ 16\cdot 8= \mtt{e}(S)\cdot n(S).
\end{align*}
Finally, let us point out that Conjecture~\ref{C:relatingconjectures} holds true for this example:
\begin{align*}
 \frac{ \mtt{c}_M(S)}{c(S)} =\frac{36}{34} \approx 1.0588 \ \leq \ \frac{ \mtt{n}_M(S)}{n(S)} = \frac{12}{8} =1.5.
\end{align*}

\end{Example}

\newpage
\subsection*{Acknowledgements}  
This research is partially supported by a grant from the Louisiana Board of Regents (090ENH-21).  
We are grateful to the anonymous referee for the helpful comments which greatly improved the quality of our paper.
The first author thanks Okinawa Institute of Science and Technology TSVP program for providing wonderful working 
conditions, where the final part of this research was completed.

\subsection*{Ethics declarations/Conflict of interests}
The authors have no conflicts of interest to declare that are relevant to the content of this article.

\subsection*{Data availability statement}
Data sharing not applicable to this article as no datasets were generated or analyzed during the current study.

\bibliography{references.bib}
\bibliographystyle{plain}
\end{document}